\documentclass{article}[12pt]

\usepackage{latexsym,color}
\usepackage{graphicx}
\usepackage{mathptmx}
\usepackage{amssymb}
\usepackage{amsfonts}
\usepackage{amsmath}
\usepackage{latexsym}
\usepackage{appendix}

\def\be{\begin{eqnarray}}
\def\ee{\end{eqnarray}}
\def\b*{\begin{eqnarray*}}
\def\e*{\end{eqnarray*}}
\def\bee{\begin{equation}}
\def\eee{\end{equation}}

\newtheorem{Theorem}{Theorem}[part]
\newtheorem{Definition}{Definition}[part]
\newtheorem{Proposition}{Proposition}[part]

\newtheorem{Assumption}{Assumption}[part]
\newtheorem{Lemma}{Lemma}[part]

\newtheorem{Remark}{Remark}[part]

\makeatletter \@addtoreset{equation}{section}

\@addtoreset{Definition}{section}

\@addtoreset{Theorem}{section}

\@addtoreset{Proposition}{section}

\@addtoreset{Property}{section}

\@addtoreset{Assumption}{section}

\@addtoreset{Corollary}{section}

\@addtoreset{Lemma}{section}

\@addtoreset{Remark}{section}

\@addtoreset{Example}{section}

\def \E{\mathbb{E}}

\def \P{\mathbb{P}}

\def \R{\mathbb{R}}

\def\Ac{{\cal A}}

\def\Ic{{\cal I}}
\def\Jc{{\cal J}}

\def\Lc{{\cal L}}
\def\Mc{{\cal M}}
\def\Nc{{\cal N}}

\def\Rc{{\cal R}}

\def\Ut{{\tilde U}}

\def \Frac{\displaystyle\frac}

\def\no{\noindent}
\def\x{\times}
\def\reff#1{{\rm(\ref{#1})}}

\def\1{{\bf 1}}
\def \qed{\hbox{ }\hfill{ ${\cal t}$~\hspace{-5.1mm}~${\cal u}$   } }

\def \proof{{\noindent \bf Proof. }}
\def \ep{\hbox{ }\hfill$\Box$}

\def\ve{v^\epsilon}
\def\ep{\epsilon}
\def\eps{\epsilon}
\def\ue{u^\epsilon}
\def\vp{v^\prime}
\def\vpp{v^{\prime \prime}}

\def\b*{\begin{eqnarray*}}
\def\e*{\end{eqnarray*}}

\def\xe{x_\ep}
\def\ye{y_\ep}
\def\ze{z^\ep}

\def\Ds{\mathbf{D}_s}
\def\Dss{\mathbf{D}_{ss}}
\def\Dx{\mathbf{D}_x}
\def\Dy{\mathbf{D}_y}

\def\Dyy{\mathbf{D}_{yy}}

\def\Dsz{\mathbf{D}_{sz}}
\def\Dsy{\mathbf{D}_{sy}}

\def\cbf{{\mathbf{c}}}
\def\ybf{{\mathbf{y}}}

\title{Homogenization and asymptotics for small transaction costs }

\author{ 
H. Mete {\sc Soner}\footnote{ETH (Swiss Federal Institute of Technology),
Zurich, and Swiss Finance Institute,
hmsoner@ethz.ch. Research partly supported by the
European Research Council under the grant 228053-FiRM
and by the ETH Foundation.}
      \and Nizar {\sc Touzi}\footnote{CMAP, Ecole Polytechnique Paris, nizar.touzi@polytechnique.edu.
      Research supported by the Chair {\it Financial Risks} of the {\it Risk Foundation} sponsored by Soci\'et\'e
             G\'en\'erale, the Chair {\it Derivatives of the Future} sponsored by the {F\'ed\'eration Bancaire Fran\c{c}aise}, and
             the Chair {\it Finance and Sustainable Development} sponsored by EDF and Calyon. }}
             
\date{\today}
 
\begin{document}

\maketitle

\begin{abstract}

We consider the classical Merton problem of lifetime 
consumption-portfolio optimization problem with
small proportional transaction costs.
The first order term in the asymptotic expansion
is explicitly calculated through a singular ergodic control
problem which can be solved in closed form 
in the one-dimensional case.  Unlike the existing literature, 
we consider a general utility function
and general  dynamics for the underlying assets.
Our arguments are based on ideas from the homogenization theory and 
use the convergence tools from the theory of viscosity solutions. 
The multidimensional case is studied in our accompanying paper 
\cite{pst} using the same approach.

\vspace{10mm}


\vspace{5mm}

\end{abstract}
\newpage

\section{Introduction}

The problem of investment and consumption in a market with transaction costs was first studied by Magill \& Constantinides  \cite{mc1976} and later by Constantinides \cite{con;86}. Since then,
starting with the classical paper of Davis \& Norman \cite{da;no;90} an impressive understanding of
this problem has been achieved. 
In these papers and in \cite{dl1991,ss} the dynamic programming approach 
in one space dimension has been developed.  
The problem of proportional transaction
costs is a special case of a singular stochastic control
problem in which the state process can have controlled discontinuities.
The related partial differential equation for this class of optimal control
problems is a quasi-variational inequality which contains a
gradient constraint.  Technically, the multi-dimensional
setting presents intriguing free boundary problems
and the only regularity result to date are 
\cite{ss1} and \cite{ss2}.
For the financial
problem, we refer to
the recent book by Kabanov \& Safarian \cite{ks2009}. 
It provides an excellent exposition 
to the later developments and the solutions in multi-dimensions.
 
It is well known that in practice the proportional transaction costs are small 
and in the limiting case of zero costs, one recovers the classical problem of 
Merton \cite{merton}. Then, a natural  approach to simplify the problem 
is to obtain an asymptotic expansion in terms of the small transaction costs.
This was initiated in the pioneering paper of Constantinides \cite{con;86}.  
The first proof in this direction was obtained in the appendix of \cite{ss}.
Later several rigorous results \cite{b2011,gms2011,js2004,r2004} and formal 
asymptotic results \cite{am2004,go2010,wh;wi;95} have been obtained. 
The rigorous results have been restricted to one space dimensions with 
the exception of the recent manuscript by Bichuch and Shreve \cite{bs2011}. 

In this and its accompanying paper \cite{pst}, we consider this classical 
problem of small proportional transaction costs and develop a unified 
approach to the problem of asymptotic analysis.
We also relate the first order asymptotic expansion in $\eps$ to an ergodic 
singular control problem.

Although our formal derivation in Section \ref{s.formal}
and the analysis of \cite{pst} are multi-dimensional,
to simplify the presentation,
in this introduction we restrict ourselves to a single risky asset
with a price process $\{S_t,t\ge 0\}$.
We assume $S_t$  is given by a time homogeneous 
stochastic differential equation together with $S_0=s$
and volatility function $\sigma(\cdot)$.
For an initial capital $z$, the value function of the 
Merton infinite horizon optimal consumption-portfolio problem 
(with zero-transaction costs) is denoted by $v(s,z)$.
On the other hand,
the value function for the problem with transaction costs
is a function of $s$ and the pair $(x,y)$ representing the wealth 
in the saving and in the stock  accounts,
respectively.  Then, the 
total wealth is simply given by $z=x+y$.
For a small proportional 
transaction cost $\ep^3>0$, we let 
$v^\ep(s,x,y)$ be the maximum expected discounted utility
from consumption.  It is clear that $v^\ep(s,x,y)$ converges to
$v(s,x+y)$ as $\ep$ tends to zero.  Our main analytical objective
is to obtain an expansion for $v^\ep$ in the small parameter
$\ep$.

To achieve such an expansion,
we assume that $v$ is smooth and let 
 \be\label{eta}
 \eta(s,z)
 &:=& 
 - \frac{v_z(s,z)}{v_{zz}(s,z)}
 \ee
be the corresponding risk tolerance. The solution of the Merton problem 
also provides us an optimal feedback portfolio strategy $\ybf(s,z)$ and an 
optimal feedback consumption function
$\cbf(s,z)$.   Then, the first term in the asymptotic expansion
is given through an ergodic singular control problem defined for every fixed point $(s,z)$ by 
 $$
 \bar{a}(s,z)
 :=
 \inf_{M} J(s,z,M),
$$
where $M$ is a control process of bounded variation with variation norm $\|M\|$, 
 $$
 J(s,z,M)
 := \limsup_{T\to\infty}\frac{1}{T}
    \E\left[\int_0^T \frac{|\sigma(s)\xi_t|^2}{2}
                     +\|M\|_T
      \right],
$$
and the controlled process $\xi$ satisfies the dynamics driven by a Brownian motion 
$B$, and parameterized by the fixed data $(s,z)$:
\b*
d \xi_t
=\alpha(s,z)
 d B_t
 +dM_t
&\mbox{where}&
\alpha:=\sigma[\ybf(1-\ybf_z)-s\ybf_s].
\e*
The above problem is defined
more generally in Remark \ref{r.ergodic}
and solved explicitly in the subsection \ref{s.w1d} below in terms of the 
zero-transaction cost value function $v$.

Let  $\{\hat Z^{s,z}_t,t\ge 0\}$ be the optimal wealth 
process using the feedback strategies $\ybf, \cbf$, and starting 
from the initial conditions $S_0=s$
and $\hat Z^{s,z}_0=z$.
Our main result 
is on the convergence of the function
 $$
 \bar{u}^\ep(x,y)
 :=
 \frac{v(s,x+y)-v^\eps(s,x,y)}{\eps^2}.
 $$

\no {\bf Main Theorem.}\quad {\it Let
$\bar{a}$ be as above and set $a:=\eta v_z\bar{a}$.
Then, as $\eps$ tends to zero,
 \be\label{e.defu}
 \bar{u}^\ep(x,y)
 \rightarrow
 u(s,z)
 :=
 \E\left[\int_0^\infty e^{-\beta t}a(S_t,\hat Z^{s,z}_t)dt\right],
 ~\mbox{locally uniformly.}
 \ee
}

Naturally, the above result requires assumptions and we refer the reader to 
Theorem \ref{t.main} for a precise statement.
Moreover, the definition and the convergence of $\ue$ is equivalent
to the expansion
\begin{equation}
\label{e.mainexp}
\ve(s,x,y) = v(z) - \ep^2 u(z) + \circ(\ep^2),
\end{equation}
where as before $z=x+y$ and
$\circ(\ep^k)$ is any function such that $\circ(\ep^k)/\ep^k$
converges to zero locally uniformly.

A  {\em formal multi-dimensional derivation} of this result is provided
in Section \ref{s.formal}.
Our approach is similar to all formal studies starting from 
the initial paper by Whalley \& Willmont 
\cite{wh;wi;95}. These formal calculations also  
provide the connection with another important class 
of asymptotic problems, namely homogenization. Indeed,
the dynamic programming equation of the ergodic problem described above
is the {\em corrector (or cell) equation} in the homogenization terminology. This identification
allows us to construct a rigorous proof similar to the ones in homogenization.
These assertions are formulated into a formal theorem at the end of Section \ref{s.formal}. 
The analysis of Section \ref{s.formal}
is very general and can easily extend to other similar problems. Moreover, the above ergodic
problem is a singular one and we show in \cite{pst} that its continuation region also describes the asymptotic shape of the no-trade region in the transaction cost problem.

The connection between homogenization and asymptotic problems in finance has already played an important role in several other problems. Fouque, Papanicolaou \& Sircar \cite{fps} use this approach for stochastic volatility models. We refer to the recent book \cite{fpss} for information on this problem and also extensions to multi dimensions. In the stochastic volatility context the homogenizing (or the so-called fast variable) is the volatility and is given exogenously. Indeed, for homogenization problems, the fast variable is almost always given. In the transaction cost problem, however, this is not the case and the main difficulty is to identify the ``fast'' variable. A similar difficulty is also apparent in a problem  with an illiquid 
financial market which becomes asymptotically liquid. 
The expansions for that problem was obtained in \cite{pst1}. We use their techniques in an essential way.

The later sections of the paper are concerned with the rigorous proof. The main technique is the viscosity approach of Evans to homogenization \cite{ev1,ev2}.
This powerful method combined with the relaxed limits of Barles \& Perthame \cite{bp} provides the necessary tools. As well known, this approach has the advantage of using only a simple  local $L^\infty$ bound which is described in Section \ref{s.assume}. In addition to \cite{bp,ev1,ev2}, the rigorous proof
utilizes several other techniques from the theory of viscosity solutions developed in the papers \cite{bp,fs89,fso,lsst,sszj,son93} for asymptotic analysis. 

For the rigorous proof, we concentrate on the simpler one dimensional setting. This simpler setting allows us to highlight the technique with the least possible technicalities. The more general multi-dimensional problem is considered in \cite{pst}.

The paper is organized as follows. The problem is introduced in the next section and the approach is formally introduced in Section \ref{s.formal}. In one dimension, the corrector equation is solved in the next section. We state the general assumptions in Section \ref{s.assume} and prove the convergence result in Section \ref{s.convergence}. In Section \ref{s.verify} we discuss the assumptions. Finally
a short summary for the power utility is given in the final Section.

\section{The general setting}
\label{s.setup}

The structure we adopt is the one developed and studied in the recent book
by Kabanov \& Safarian \cite{ks2009}. We briefly recall it here.

We assume a financial market consisting of a non-risky asset $S^0$ 
and $d$ risky assets with price process $\{S_t=(S^1_t, \ldots,S^d_t),t\ge 0\}$ 
given by the stochastic differential equations,
 $$
 \frac{dS^0_t}{S^0_t}
 =
 r(S_t)dt,~~
 \frac{dS^i_t}{S^i_t}
 =
 \mu^i(S_t) dt + \sum_{j=1}^d \sigma^{i,j}(S_t) dW^j_t,
 ~~1\le i\le d,
 $$
where $r:\R^d\to \R_+$ is the instantaneous interest rate and $\mu:\R^d \to \R^d$, 
$\sigma:\R^d\to\Mc_{d}(\R)$ are the coefficients of instantaneous mean return and volatility.
We use the notation $\Mc_{d}(\R)$ to denote 
$d\times d$ matrices with real entries. 
The standing assumption on the coefficients
 \b*
 r, \mu, \sigma &\mbox{are bounded and Lipschitz, and}&
 (\sigma \sigma^{T})^{-1}~~\mbox{is bounded,}
 \e*
will be in force throughout the paper (although not recalled in our statements). 
In particular, the above stochastic differential equation has a unique strong solution.

The portfolio of an investor is represented by the 
dollar value $X$ invested in the non-risky asset
 and the vector process $Y=(Y^1,\ldots,Y^d)$ of 
the value of the positions in each risky asset.
 The portfolio position is allowed to change in 
 continuous-time by  transfers from any asset to 
 any other one. 
 However, such transfers are subject to 
 {\em{proportional transaction costs}}. 

We continue by describing the portfolio rebalancing 
in the present setting.
For all $i,j=0,\ldots,d$, let $L^{i,j}_t$ be
the total amount of transfers (in dollars)
from the $i$-th to the $j$-th asset cumulated up to time $t$. 
Naturally, the processes $\{L^{i,j}_t,t\ge 0\}$ are defined
 as c\`ad-l\`ag, nondecreasing, adapted processes with $L_{0^-}=0$
 and $L^{i,i}\equiv 0$. 
 The proportional transaction cost induced by a 
 transfer from the $i$-th to the $j$-th stock is 
 given by $\eps^3\lambda^{i,j}$ where $\eps>0$ is a small parameter, and  
 $$
 \lambda^{i,j}\ge 0,
 ~~
 \lambda^{i,i}=0,
 ~~i,j=0,\ldots,d.
 $$
The scaling $\ep^3$ is chosen to state the expansion results simpler.
We refer the reader to the recent book of Kabanov \& Safarian  \cite{ks2009} 
for a thorough discussion of the model.
 
The {\em solvency region} $K_\eps$ is defined as the set of 
all portfolio positions which can be
 transferred into portfolio positions with nonnegative 
 entries through an appropriate portfolio rebalancing.
We use the notation $\ell=(\ell^{i,j})_{i,j=0,\ldots d} $ to denote 
this appropriate 
instantaneous transfers of size $\ell^{i,j}$.
We directly compute 
that the induced change in each entry, after subtracting the 
corresponding transaction costs is given by the linear operator
$\mathbf{R}:\Mc_{d+1}(\R_+)\to\R^{d+1}$,
 \b*
 \mathbf{R}^i(\ell)
 :=
 \sum_{j=0}^{d}
 \big(\ell^{j,i}-(1+\eps^3\lambda^{i,j})\ell^{i,j}\big),
 ~~i=0,\ldots,d,
 &\mbox{for all}&
 \ell\in\Mc_{d+1}(\R_+),
 \e* 
where $\ell^{i,j}>0$ and $\ell^{j,i}>0$ for some $i,j$ would clearly be suboptimal.
Then, $K_\ep$ is given by
 \begin{eqnarray*}
 &&
 \!\!\!\!
 K_\eps
 :=
 \Big\{(x,y)\in\R \times \R^{d}:~
       (x,y)+\mathbf{R}(\ell)
       \in \R_+^{1+d}
       ~~\mbox{for some}~~\ell\in\Mc_{d+1}(\R_+) 
 \Big\}.
 \end{eqnarray*}
For later use, we denote by $(e_0,\ldots,e_d)$ the 
canonical basis of $\R^{d+1}$ and set
$$
 \Lambda_{i,j}^\eps
 :=
 e_i-e_j+\eps^3\lambda^{i,j}\;e_i,
 \qquad  i,j=0,\ldots, d.
 $$
 
In addition to the trading activity, the investor consumes at a rate 
determined by a nonnegative progressively measurable 
process $\{c_t,t\ge 0\}$. Here $c_t$ represents the rate of 
consumption in terms of the non-risky asset $S^0$. 
Such a pair $\nu:=(c,L)$ is called a 
{\em consumption-investment strategy}. 
For any initial position $(X_{0^-},Y_{0^-})=(x,y)\in\R\x\R^{d}$, 
the portfolio position of the investor are 
given by the following state equation
 $$
 dX_t
 =
 \big(r(S_t)X_t-c_t\big)dt
 +\mathbf{R}^0(dL_t),
 ~~\mbox{and}~~
 dY^i_t
 =
  Y^i_t\;\frac{dS^i_t}{S^i_t}
 +\mathbf{R}^i(dL_t),
 ~~i=1,\ldots,d.
 $$
 The above solution depends
 on the initial condition $(x,y)$,
 the control $\nu$ and also on the initial
 condition of the stock process $s$.
Let  $(X,Y)^{\nu,s,x,y}$
 be the solution of the above equation.
 Then, a consumption-investment strategy $\nu$ is said to 
be {\em admissible} for the initial position $(s,x,y)$, if
$$
(X,Y)^{\nu,s,x,y}_t\in K_\eps,
 \qquad \forall \ t\ge 0,
 \quad
 \P-\mbox{a.s.}
 $$
 The set of admissible strategies
 is denoted by $\Theta^\ep(s,x,y)$.
For given initial
positions $S_{0}=s \in \R_+^d$, 
$X_{0^-}=x \in \R$,
$Y_{0^-}=y \in \R^d$,
the investment-consumption problem is
the following maximization problem,
 \b*
 \ve(s,x,y)
 &:=&
 \sup_{(c,L) \in \Theta^\ep(s,x,y)}\
 \E\left[\int_0^\infty\ e^{-\beta t}\ U(c_t)dt \right],
 \e*
where 
$U:(0,\infty)\mapsto \R$ is a utility function. 
We assume that $U$ is $C^2$, increasing, strictly concave, 
and we denote its convex conjugate by,
 \b*
 \Ut(\tilde c)
 &:=& 
 \sup_{c>0} \big\{U(c)-c\tilde c\big\},
 \qquad
 \tilde c\in\R.
 \e*
Then $\tilde U$ is a $C^2$ convex function.
It is well known that the value function
is a viscosity solution of
the corresponding dynamic programming equation.
In one dimension,
this is first proved in \cite{ss}.
In the above generality, we refer to 
\cite{ks2009}. 
To state the equation, we first need to introduce 
some more notations. We define a second order 
linear partial differential operator by,
 \be\label{Lc}
 \Lc
 &:=&
 \mu\cdot\left(\Ds +\Dy\right)
 +r\Dx
 +\frac12\mbox{Tr}\left[\sigma \sigma^{\rm T}\left(\Dyy+\Dss+2 \Dsy\right) \right],
 \ee
where $^{\rm{T}}$ denotes the transpose and for $i,j=1,\ldots,d$,
 \b*
 &\Dx:=x\Frac{\partial}{\partial x},
 ~~\Ds^i:=s^i\Frac{\partial}{\partial s^i},
 ~~\Dy^i:=y^i\Frac{\partial}{\partial y^i},&
 \\
 &\Dss^{i,j}:=s^is^j\Frac{\partial^2}{\partial s^i\partial s^j},
 ~~\Dyy^{i,j}:=y^iy^j\Frac{\partial^2}{\partial y^i\partial y^j},
 ~~\Dsy^{i,j}:=s^iy^j\Frac{\partial^2}{\partial s^i\partial y^j},&
 \e*
  $\Ds=(\Ds^i)_{1\le i\le d}$, $\Dy=(\Dy^i)_{1\le i\le d}$, 
  $\Dyy:=(\Dyy^{i,j})_{1\le i,j\le d}$, $\Dss:=(\Dss^{i,j})_{1\le i,j\le d}$, 
  $\Dsy:=(\Dsy^{i,j})_{1\le i,j\le d}$.
Moreover, for a smooth scalar function $(s,x,y)\in \R^d_+\times\R\times\R^d
\mapsto \varphi(x,y)$, we set
 \b*
 \varphi_{x}:=\frac{\partial\varphi}{\partial x}\;\in\R,
\qquad
 \varphi_{y}:=\frac{\partial\varphi}{\partial y}\;\in\R^d.
 \e*

\begin{Theorem}
\label{t.dpp}
Assume that the value function $\ve$ is locally bounded. 
Then, $\ve$ is a viscosity solution of the dynamic programming equation in 
$\R_+^d\x K_\ep$,
 \be \label{e.dpp}
 \min_{0\le i,j\le d}
 \left\{\ \beta \ve - \Lc \ve - \Ut(\ve_{x})\ , 
        \ \Lambda_{i,j}^\eps\cdot(\ve_x,\ve_y) \ 
 \right\} 
 =
 0.
 \ee
Moreover, $\ve$ is concave in $(x,y)$ and converges to the Merton value function 
$v:=v^0$, as $\eps>0$ tends to zero.
\end{Theorem}

Under further conditions the uniqueness in the above statement is proved in 
\cite{ks2009}. However, this is not needed in our subsequent analysis. 

\subsection{Merton Problem}

The limiting case of $\eps=0$ 
corresponds to the 
classical Merton portfolio-investment problem
 in a frictionless financial market. 
In this  limit, since the transfers
 from one asset to the other are costless, 
 the value of the portfolio can be measured 
 in terms of the nonrisky asset $S^0$. 
 We then denote by $Z:=X+Y^1+\ldots+Y^d$
  the total wealth obtained by 
the aggregation of the positions on all assets. 
In the present setting, we denote by $\theta^i:=Y^i$ and 
$\theta:=(\theta^1,\ldots,\theta^d)$ the vector 
process representing the positions on the risky assets. 
The wealth equation for the Merton problem is 
then given by
 \be
 dZ_t 
 =
 \big(r(S_t)Z_t-c_t\big)dt
 +\sum_{i=1}^d\ \theta^i_t \;\Big(\frac{dS^i_t}{S^i_t}-r(S_t)dt
                             \Big).
 \ee
An admissible consumption-investment strategy is
 now defined as a pair $(c,\theta)$ of progressively 
 measurable processes with values in $\R_+$ and $\R^d$, 
 respectively, and such that the corresponding 
 wealth process is well-defined and 
 almost surely non-negative for all times. The set of all
  admissible consumption-investment strategies is denoted by $\Theta(s,z)$.
 
The Merton optimal consumption-investment 
problem is defined by
 \b*
 v(s,z)
 := 
 \sup_{(c,\theta) \in \Theta(s,z)}\
 \E\left[\int_0^\infty\ e^{-\beta t}\ U(c_t)dt \right],
 &s\in \R_+^d,&
 z\ge 0.
 \e*
Throughout this paper, we assume that the Merton value function 
$v$ is strictly concave in $z$ and is a classical solution of 
the dynamic programming equation,
$$
 \beta v - rz v_z - \Lc^0 v - \Ut(v_z)
 -\sup_{\theta\in\R^d}\Big\{\theta\cdot\big((\mu-r\1_d) v_z
                                             +\sigma \sigma^{\rm T}\Dsz v
                                       \big)
                              +\frac12|\sigma^{\rm T}\theta|^2 v_{zz}
                        \Big\}
 =
 0,
$$
where $\1_d:=(1,\ldots,1)\in \R^d$, $\Dsz:=\frac{\partial}{\partial z}\Ds$, and
 \be\label{Lc0}
 \Lc^0
 &:=&
 \mu\cdot\Ds
 +\frac12\mbox{Tr}\big[\sigma \sigma^{\rm{T}} \Dss\big].
 \ee
The optimal controls are smooth functions $\cbf(s,z)$ and $\ybf(s,z)$ 
obtained by as the maximizers of the Hamiltonian. Hence, 
 \begin{equation}
 \label{e.merton1}
 0 =\beta v - \Lc^0v - \Ut(v_z)
 -rzv_z 
 -\ybf\cdot(\mu- r \1_d) v_z
 -\sigma \sigma^{\rm T} \ybf \cdot\Dsz v
 -\frac12|\sigma^{\rm T}\ybf|^2 v_{zz},
 \end{equation}
the optimal consumption rate is given by,
 \be\label{hatc}
 \cbf(s,z)
 :=
 -\Ut^{\prime}\big(v_z(s,z)\big)
 =\big(U^{\prime}\big)^{-1}\big(v_z(s,z)\big)
 \ \mbox{for}\ \ \
 s\in\R_+^d,\ \ z\ge 0,
 \ee
and the optimal investment strategy  $\ybf$
is obtained by solving the finite-dimensional maximization problem,
 $$
 \max_{\theta\in\R^d}
 \Big\{\frac12 |\sigma^{\rm T}\theta|^2 v_{zz}
       +\theta\cdot\big((\mu-r\1_d) v_z+\sigma \sigma^{\rm T}\Dsz v\big)
 \Big\}.
 $$
Since $v$ is strictly concave, the Merton 
optimal investment strategy $\ybf(s,z)$ satisfies
 \be\label{hattheta}
 -v_{zz}(s,z)\;\sigma\sigma^{\rm T}(s)\ybf(s,z)
 =
 (\mu-r\1_d)(s) v_z(s,z)+\sigma\sigma^{\rm T}(s) \Dsz v(s,z).
 \ee

\section{Formal Asymptotics}
\label{s.formal}

In this section, we provide the  formal derivation of the expansion
in any space dimensions. In the subsequent sections, we prove this expansion 
rigorously for the one dimensional case. Convergence proof in higher 
dimensions is carried out in a forthcoming paper \cite{pst}.  
In the sequel we use the standard notation $O(\ep^k)$ to denote any function which 
is less than a locally bounded function times $\ep^k$ and $\circ(\ep^k)$
 is a function such that $\circ(\ep^k)/\ep^k$ converges to zero locally uniformly.

Based on previous results \cite{wh;wi;95,am2004,go2010,js2004,ss}, 
we postulate the following expansion,
 \begin{equation}
 \label{e.expansion}
 \ve(s,x,y) = v(s,z) - \ep^2 u(s,z) - \ep^4 w(s,z,\xi) +\circ(\ep^2),
 \end{equation}
where
$(z,\xi)=(z,\xi_\ep)$ is 
a transformation of 
$(x,y)\in K_\eps$ given by
 \begin{equation}
 \nonumber
 z= x+y^1+\ldots+y^d,
 \qquad
 \xi^i:=\xi^i_\ep(x,y) = \frac{y^i- \ybf^i(s,z)}{\ep},
 \qquad i=1,\ldots, d,
 \end{equation}
$\ybf=\big(\ybf^1,\ldots,\ybf^d\big)$ is the Merton optimal investment 
strategy of \reff{hattheta}. In the postulated expansion 
\reff{e.expansion}, we have also introduced two functions
 \b*
 u:\R_+^d\times\R_+ \mapsto \R,
 &\mbox{and}&
 w:\R_+^d\times\R_+ \x \R^d \mapsto \R.
 \e*
The main goal of this section is
to formally derive equations for these two functions. 
A rigorous proof will be also provided in the 
subsequent sections and the precise  
statement for this expansion is 
stated in  Section \ref{s.convergence}.

Notice that the above expansion is assumed to hold up to $\ep^2$, i.e.
the $\circ(\ep^2)$ term. Therefore, the reason for having a higher term 
like $\ep^4w(z,\xi)$ explicitly in the expansion may not be clear. 
However, this term contains the fast variable $\xi$ and its second 
derivative is of order $\ep^2$, which will then contribute to the asymptotics 
since $v^\eps$ solves a second order PDE. This follows the intuition introduced 
in the pioneering work of Papanicolaou and 
Varadhan \cite{pv79} in the theory of homogenization.

Since $(x,y)\in K_\eps\mapsto(z,\xi)\in\R_+\times \R^d$ 
is a one-to-one change of variables, in the sequel for any function 
$f$ of $(s,x,y)$ we use the convention,
 \be
 \label{e.hat}
 \hat{f}(s,z,\xi):= f\big(s,z-\ep \xi-\ybf(s,z),\ep \xi+\ybf(s,z)\big).
 \ee
The new variable $\xi$ is the ``fast'' 
variable and in the 
limit it homogenizes to yield the
convergence of $\hat{v}^\ep(s,z,\xi)$
 to the Merton function $v(s,z)$ which 
 depends only on the $(s,z)$-variables.  
 This is the main formal connection of 
 this problem to the theory of homogenization.
 This variable was also used
 centrally by Goodman \& Ostrov \cite{go2010}.
 Indeed, their asymptotic results
 use the properties of the stochastic 
 equation satisfied by $\ep\xi^\ep(X_t,Y_t)$. 

First we directly differentiate the expansion \reff{e.expansion}
and compute the terms appearing in \reff{e.dpp} in term of $u$ and $w$.
The directional derivatives are given by,
$$
 \Lambda_{i,j}^\eps\cdot (v^\ep_{ x},v^\ep_{ y}) =
 -\eps^4(e_i-e_j)\cdot (w_x(s,z,\xi),w_y(s,z,\xi))
 +\eps^3\lambda^{i,j}v_z
 +O(\eps^4).
$$
We directly calculate that,
 \begin{equation}
 \label{wx}
(w_x,w_y)(s,z,\xi)
 =\Big(w_z-\frac{1}{\eps}\ybf_z\cdot w_\xi\Big)\1_{d+1}
 +\frac{1}{\eps}
  \left(0,w_\xi
  \right).
 \end{equation}
To simplify the notation, we introduce
\be
\label{e.dhat}
\hat{D}_\xi w(s,z,\xi):= (0, {D}_\xi w(s,z,\xi)) \in \R^{d+1}.
\ee
Then,
 \be\label{e.dpm}
 \Lambda_{i,j}^\eps\cdot(v^\ep_{ x},v^\ep_{ y}) 
 =
 \ep^3\big( \lambda^{i,j}v_z
 +(e_j-e_i)\cdot \hat{D} w) 
 + O(\ep^4).
 \ee
The elliptic equation in \reff{e.dpp} requires a longer 
calculation and we will later use the Merton identities  
\reff{e.merton1}, \reff{hatc} and \reff{hattheta}.
Firstly, by \reff{e.merton1},
 \b*
 I^\ep
 &:=& 
 \beta\ve - \Lc \ve - \Ut(\ve_x)
 \\ 
 &=& 
 (\ybf-y)\cdot \big[(\mu-r \1_d) v_z+\sigma  \sigma^{\rm T} \Dsz v\big] 
 +\frac12 \big(| \sigma^{\rm T}\ybf|^2-| \sigma^{\rm T} y|^2
          \big)v_{zz}
 \\ 
 && + \Big(\Ut(v_z) 
           - \Ut\big(v_z +\ep^2 u_z +O(\ep^3)\big)\Big)
 \\
 && - \ep^2 \Big(\beta u - {\Lc} u 
            \Big) 
    +\frac{\eps^4}{2}  {\mbox{Tr}}[\sigma  \sigma^{\rm T}\Dyy w]
    +O(\ep^3).
 \e*
We use Taylor expansions on the terms 
involving $\Ut$ and \reff{hatc}-\reff{hattheta} in the first line,
to arrive at
 \be
 I^\eps
 &=&
 \Big(-\sigma^{\rm T}(\ybf-y)\cdot \sigma^{\rm T}\ybf
 +\frac12 \big(|\sigma^{\rm T}\ybf|^2-|\sigma^{\rm T} y|^2
          \big)\Big)v_{zz}
 \nonumber \\ 
 && - \ep^2 \Big(\beta u - {\Lc} u + \hat c u_z 
            \Big) 
    +\frac{\eps^4}{2} 
    {\mbox{Tr}}[\sigma \sigma^{\rm T} (\Dyy+\Dss+\Dsy)w]
    +O(\ep^3)
 \nonumber\\
 &=&
 -\frac12{|\sigma^{\rm T}(\ybf-y)|^2}  v_{zz}
 - \ep^2 \Big(\beta u - {\Lc} u + \hat c u_z 
         \Big) 
 +\frac{\eps^4}{2}  
 {\mbox{Tr}}[\sigma \sigma^{\rm T}(\Dyy+\Dss+\Dsy)w]
 +O(\ep^3)
 \nonumber\\
 &=&
 \ep^2 \Big( -\frac12{|\sigma^{\rm T}\xi|^2} v_{zz}
             - \beta u + {\Lc} u - \hat c u_z 
       \Big) 
 +\frac{\eps^4}{2} 
 {\mbox{Tr}}[\sigma \sigma^{\rm T}(\Dyy+\Dss+\Dsy)w]
 +O(\ep^3).
 \label{e.ie0}
 \ee
Finally, from \reff{wx}, we see that
$$
 \partial_yw 
 =
 w_z\1_d 
 +\frac{1}{\eps}\big(I_d-\1_d\ybf^{\rm T}_z\big) w_\xi.
$$
Therefore,
 $$
 \partial_{yy}w
 =
 \Big(w_{z z}
 -\frac{1}{\ep}(\ybf_{zz}\cdot w_\xi+\ybf_z\cdot w_{z\xi})\Big)
               \1_d \1_d^{\rm T}
 +\frac{1}{\ep}\big(w_{z\xi}\1_d^{\rm T}+\1_dw_{z\xi}^{\rm T}\big)
 +\frac{1}{\ep^2}\big(I_d-\1_d\ybf^{\rm T}_z\big)
                 w_{\xi \xi}
                 \big(I_d-\ybf_z\1_d^{\rm T}\big).
 $$
We substitute this in \reff{e.ie0} and use the fact that 
$y=\ybf+O(\eps)$. This yields,
 \be
 I^\eps
 =
 \ep^2 \Big( -\frac12{|\sigma^{\rm T}\xi|^2} v_{zz}
             +\frac12 {\mbox{Tr}}\big[\alpha \alpha^{\rm T} w_{\xi\xi}\big]
             -\Ac u
       \Big) 
 +O(\ep^3),
 \label{e.ie1}
 \ee
where $\alpha(s,z)$ is given by
\begin{equation}\label{alpha}
 \alpha(s,z)
 =
 \big\{\big(I_d-\ybf_z\1_d^{\rm T}\big) \mbox{diag}[\ybf]
       -\ybf_s^{\rm T}\mbox{diag}[s]
 \big\}(s,z)
 \sigma(s) ,
\end{equation}
$\mbox{diag}[\ybf]$ denotes the diagonal matrix 
with $i$-th diagonal entry $\ybf^i$, and
 \begin{equation}
\label{Ac}
 \Ac u 
 = \beta u 
 -\Lc^0 u
 - \big(rz + \ybf\cdot(\mu-r\1_d) - \cbf\big) 
   u_z
 -\frac12 |\sigma^{\rm T} \ybf|^2 \,u_{zz}
 -\sigma\sigma^{\rm T}\ybf\cdot\Dsz u.
 \end{equation}
Recall that $\Lc^0$ is the
infinitesimal generator of the stock price process.
Observe that the above operator is the infinitesimal generator of the pair process 
$(S,\hat{Z})$ where $\hat{Z}$ is the optimal wealth process in 
the Merton zero-transaction cost problem corresponding to 
the optimal feedback controls $(\cbf,\ybf)$. In particular, 
the dynamic programming equation \reff{e.merton1}
for the Merton problem  may be expressed as,
 \begin{equation}
 \label{e.av}
 \Ac v(s,z) = U(\cbf(s,z)).
 \end{equation}
We have now obtained expressions for all the terms in the dynamic 
programming equation \reff{e.dpp}. We substitute \reff{e.dpm} and 
\reff{e.ie1} into \reff{e.dpp}. Notice that since
$\epsilon>0$, for any $A, B$, $\max\{\eps^2 A,\eps^3 B\}=0$ 
is equivalent to $\max\{A,B\}=0$. 
Hence, $w$ and $u$ satisfy,
 \b*
 \max_{0\le i,j\le d}
 & \max & \Big\{ \frac12{\big|\sigma^{\rm T}(s)\xi\big|^2} v_{zz}(s,z)  
              -\frac12 \mbox{Tr}\big[\alpha\alpha^{\rm T}(s,z)
                                      w_{\xi\xi}(s,z,\xi)
                                \big] 
              + a(s,z)\ ,
 \\
 &&
 \hspace{40mm}
 -\lambda^{i,j}v_z(s,z) + (e_i-e_j)\cdot \hat{D}_{\xi}w(z,\xi)
 \Big\}
 \;=\;
 0.~~
\e*
where  $\hat{D}_{\xi}=(0,D_\xi w)$ is as in \reff{e.dhat}
and  $a$ is given  by,
 \b*
 a(s,z):=  \Ac u (s,z),
 &s\in \R_+^d,&
 z>0.
 \e*
In the first equation above, the pair $(s,z)$ is simply a parameter
and the independent variable is $\xi$. 
Also the value of the function $w(s,z,0)$ is irrelevant in 
\reff{e.expansion} as it only contributes to the $\ep^4$ term.  
Therefore, to obtain a unique $w$, we set its value at the origin to zero. 
We continue by presenting these equations in a form
that is compatible with the power case.  So we first
divide the above equation by $v_z$ and then introduce the 
new variable 
 \b*
 \rho &=& \xi/\eta(s,z),
 \e* 
where $\eta$ is the risk tolerance 
coefficient defined by \reff{eta}.
We also set
$$
\bar{w}(s,z,\rho):= \frac{w(s,z,\eta(s,z) \rho)}{
\eta(s,z) v_z(s,z)},
\quad
\bar{a}(s,z):= \frac{a(s,z)}{
\eta(s,z) v_z(s,z)},
\quad
\bar{\alpha}(s,z):= \frac{\alpha(s,z)}{\eta(s,z)}.
$$
Then, the {{corrector equations}} in this context is the 
following pair of equations.

\vspace{10pt}
\begin{Definition}[Corrector Equations]
{\rm{For a given point  $(s,z) \in \R_+^d \x\ R_+$,}}
the first corrector equation {\rm{is  for the unknown pair}} $(\bar{a}(s,z),
 \bar{w}(s,z,\cdot)) \in \R \times C^2(\R^d)$,
 \begin{eqnarray}
 \label{e.cw}
 &&
 \max_{0\le i,j\le d}
 \max \Big\{ -\frac{|\sigma^{\rm T}(s)\rho|^2}{2}
 -\frac12 \mbox{Tr}\big[\bar{\alpha}\bar{\alpha}^{\rm T}(s,z) 
 \bar{w}_{\rho \rho}(s,z,\rho)\big] 
              + \bar{a}(s,z)\ ,
 \\
\nonumber
 &&
\hspace{40mm}
 -\lambda^{i,j} + (e_i-e_j)\cdot \hat{D}_{\rho}\bar{w}(s,z,\rho)
 \Big\}
 \;=\;
 0,
\quad \forall \ \rho \in \R^d,
 \end{eqnarray}
  {\rm{together with the normalization $\bar{w}(s,z,0)=0$. 
 
The}} second corrector equation  {\rm{uses the constant term 
$\bar a(s,z)$ from the first corrector equation 
and it is a simple linear equation for the function $u :\R_+^d\x\R^+ 
 \mapsto \R^1$,
 \begin{equation}\label{e.au}
 \Ac u(s,z) = a(s,z) = v_z(s,z) \eta(s,z) \bar{a}(s,z), 
 \qquad
\forall \  s\in\R_+^d,~z \in \R^+.
 \end{equation}
We say that the pair $(u,w)$ is {\em{the solution of the corrector equations}} 
for a given utility function or equivalently for a given Merton value function.}}
\qed
\end{Definition}

We summarize our formal calculations in the following.  
\vspace{10pt}

\no {\bf{Formal Expansion Theorem.}}
{\em{The value function has the expansion 
\reff{e.expansion} where  $(u,w)$ is the unique solution
 of the corrector equations.
}}

\begin{Remark}{\rm The function $u$ introduced in
 \reff{e.defu} is a solution of the second corrector equation 
 \reff{e.au}, provided that it is finite. Then, assuming that 
 uniqueness holds for the linear PDE \reff{e.au} in a 
 convenient class, it follows that $u$ is given by the 
 stochastic representation \reff{e.defu}.
}
\end{Remark}
\begin{Remark}{\rm
Usually a second order equation like \reff{e.au}
in $(0,\infty)$  needs
to be completed by a boundary condition at the origin.
However, as we have already remarked, the operator $\Ac$
is the infinitesimal generator of the optimal wealth process
in the Merton problem. Then, under the Inada conditions 
satisfied by the utility function $U$, we expect that
this process does not reach the origin.  Hence, we only 
need appropriate growth conditions
near the origin and at infinity to ensure uniqueness.
}\qed
\end{Remark}

\begin{Remark}
\label{r.ergodic}{\rm{
The first corrector equation has the following stochastic representation
as the dynamic programming equation of an ergodic control problem.
For this representation we fix $(s,z)$ and
let $\{M^{i,j}_t,t\ge 0\}$ be non-decreasing control processes, 
for each $i,j=0,\ldots,d$. 
Let $\rho$ be the controlled process defined by,
 \b*
 \rho^i_t
 &=& 
 \rho^i_0
 +\sum_{j=1}^d \bar{\alpha}^{i,j}(s,z) B^j_t
 + \sum_{j=0}^d \big(M^{j,i}_t - M^{i,j}_t\big),
 \e*
for some arbitrary initial condition $\rho_0$
and a $d$ dimensional 
standard Brownian motion $B$.
Then, the ergodic control problem is
$$
\bar{a}(s,z):= \inf_{M} \ J(s,z,M),
$$
where
 $$
 J(s,z,M)
 := 
 \limsup_{T\to \infty}\ 
 \frac{1}{T} \ 
 \E\Big[ \frac12 \int_0^T \big|\sigma^{\rm T}(s)\rho_t\big|^2dt
        +\sum_{i,j=0}^d \lambda^{i,j} M^{i,j}_T 
   \Big].
 $$ 
In the scalar case,
this problem is closely related
to the classical finite fuel problem introduced by 
Benes, Shepp \& Withenhaussen \cite{bsw}.
We refer
to the paper by Menaldi, Robin and Taksar \cite{mrt}
for the present multidimensional setting.
 
The function $\bar{w}$ is the so-called
potential function in ergodic control.  We refer the reader to
the book and the manuscript of Borkar \cite{bor,bor1} 
for information on the dynamic programming approach for
the ergodic control 
problems.}}
\qed
\end{Remark}

\begin{Remark}
\label{r.formal}
{\rm
The  calculation leading to \reff{e.ie1} is used several times in the paper. 
Therefore, for future reference, we summarize it once again.
Let $v$, $z$ and $\xi$ be as above. For {\em{any}} smooth functions
 $$
 \phi:\R_+^d\x\R_+\mapsto \R,\quad
 \varpi:\R_+^d\x\R_+\x\R^d \mapsto \R,
 $$
and $\ep \in (0,1]$ set
 $$
 \Psi^\ep(s,x,y):= v(s,z) -\ep^2 \phi(s,z) - \ep^4 \varpi(s,z,\xi).
 $$
In the above calculations, we obtained an expansion 
for the second order nonlinear operator
 \be
 \nonumber
 \Jc(\Psi^\ep)
 &:=& 
 \beta\Psi^\ep 
      - \Lc \Psi^\ep 
      - \Ut(\Psi^\ep_{x})
 \\
 \label{e.rem}
 &=&
 \ep^2 \Big(-\frac{v_{zz}}{2} |\sigma^{\rm T}\xi|^2  
            +\frac12 {\mbox{Tr}}\big[\alpha\alpha^{\rm T}\varpi_{\xi\xi}\big]
            - \Ac \phi
            + \Rc^\ep
       \Big),
 \ee
where $\alpha$, $\Ac$ are as before and $\Rc^\ep(s,x,y)$ is the 
remainder term. Moreover, $\Rc^\ep$ is locally bounded by a 
$\ep$ times a constant depending only on the values of 
the Merton function $v$, $\phi$ and $\varpi$. 
Indeed, a more detailed description and an estimate 
will be proved in one space dimension in Section \ref{s.convergence}. 
} \qed
\end{Remark}

\section{Corrector Equation in one dimension}
\label{s.corrector}

In this section, we solve the first corrector 
equation explicitly in the one-dimensional case. 
Then, we provide some estimates for 
the remainder introduced in Remark \ref{r.formal}.

\subsection{Closed-form solution of the first corrector equation}
\label{s.w1d}

Recall that $w= \eta v_z\bar{w}$, $a=\eta v_z \bar{a}$, and 
the solution of the corrector equations is
a pair $(\bar{w},\bar a)$ satisfying,
 \begin{equation}\label{e.1corrector-dim1}
 \max\Big\{-\frac12\sigma^2\rho^2
           -\frac12 \bar{\alpha}^2 \bar w_{\rho \rho}
           + \bar a,
           -\lambda^{1,0}+\bar w_\rho,
           -\lambda^{0,1}-\bar w_\rho
     \Big\}
 =
 0,
 ~~\bar w(s,z,0)=0,
 \end{equation}
where $\bar{\alpha}= \alpha/\eta$ and
$\alpha(s,z)$ is given in \reff{alpha}.
We also recall that the variables $(s,z)$ are fixed
parameters in this equation. 
Therefore, throughout this section, we suppress the dependences of 
$\sigma,\alpha$ and $\bar w$ on these variables. 

In order to compute the solution
explicitly in terms of $\eta$,
we postulate a solution of the form
 \begin{equation}
 \label{wbar1d}
 \bar w(\rho)
 =
 \left\{ \begin{array}{ll}
         k_4\rho^4+k_2\rho^2+k_1\rho,& \rho_1\le\rho\le\rho_0,
         \\
         \bar w(\rho_1)-\lambda^{0,1}(\rho-\rho_1),& \rho\le\rho_1,
         \\
         \bar w(\rho_0)+\lambda^{1,0}(\rho-\rho_0),
        & \rho\ge\rho_0.
         \end{array}
 \right.
 \end{equation}
We first determine $k_4$ and $k_2$ by imposing that 
the fourth order polynomial solves the second 
order equation in $(\rho_0,\rho_1)$. 
A direct calculation yields,
 \be
\nonumber
 k_4=\frac{-\sigma^2}{12\bar\alpha^2}
 &\mbox{and}&
 k_2= \frac{\bar a}{\bar\alpha^2}.
 \ee
We now impose the smooth pasting condition, namely
assume that  $\bar w$ is $C^2$ 
at the points $\rho_0$ and $\rho_1$.
Then, the continuity of the second derivatives yield,
 \be
\label{e.xi}
 \rho_0^2=\rho_1^2=\frac{2\bar a}{\sigma^2}
 &\mbox{implying that}&
 \bar a\ge 0
 ~~\mbox{and}~~
 \rho_0=-\rho_1=\Big(\frac{2\bar a}{\sigma^2}\Big)^{1/2}.
 \ee
The continuity of the first derivatives
of $\bar w$ at the points $\rho_0$ and $\rho_1$
yield,
 \b*
 4k_4(\rho_0)^3+2k_2\rho_0+k_1 &=& -\lambda^{0,1},
 \\
 4k_4(\rho_1)^3+2k_2\rho_1+k_1 &=& \lambda^{1,0}.
 \e*
Since $\rho_0=-\rho_1$, we determine the value of 
$k_1$ by summing the two equations,
 \be
\nonumber
 k_1 &=& \frac{\lambda^{1,0}-\lambda^{0,1}}{2}.
 \ee
Finally, we obtain the value of $\bar a$ by further substituting the 
values of $k_4$, $k_2$ and $\rho_0=-\rho_1$.
The result is
 \be
\label{e.bara}
 \bar a
 =
 \frac{\sigma^2}{2}\rho_0^2
 &\mbox{and}&
 \rho_0=\Big(\frac{3\bar\alpha^2}{4\sigma^2}
                  (\lambda^{1,0}+\lambda^{0,1})
             \Big)^{1/3}.
 \ee
All coefficients of our candidate are now uniquely determined.
Moreover,  we verify that the gradient constraint
 \be
\label{bound-wxi}
 -\lambda^{1,0}\le \bar w_\rho\le \lambda^{0,1}
 \ee
holds true for all $\rho\in\R$. Hence, $\bar{w}$
constructed above is a solution of
the corrector equation. 
One may also prove that
it is the unique solution.  
However, in the subsequent analysis 
we simply use the function $\bar{w}$ defined in 
\reff{wbar1d} with the constants determined above.
Therefore, we do not study the question of uniqueness
of the corrector equation.

\begin{Remark}{\rm{
In the homothetic case with constant coefficients $r,\mu$, and $\sigma$, 
one can explicitly calculate all the functions, see Section \ref{s.homothetic}.  
Here we only report that, in that case, 
all functions are independent of the $s-$variable 
and $\rho_0, \bar{a}(z)$ are constants. 
Therefore, $a(z)$ is a positive constant times the Merton value function.}}
\qed
\end{Remark}

\begin{Remark}
\label{e.west}{\rm  
Pointwise estimates on the derivatives of 
$w$ will be used in the subsequent sections. 
So we record them here for future references. 
Indeed, by \reff{bound-wxi} 
and the fact that $w(\cdot,0)=0$,  
 \begin{equation}
\label{boundwwxi}
|w(s,z,\xi)| \le \overline{\lambda}\,v_z(s,z)|\xi|,\ \
 |w_\xi(s,z,\xi)| \le \overline{\lambda}\,v_z(s,z),\ \
 \mbox{where}\ \
 \overline{\lambda}:=\lambda^{0,1}\vee\lambda^{1,0}.
 \end{equation}
Moreover, under the smoothness assumption on $v$, 
we obtain the following pointwise estimates
 \be
 &\big(|w|+|w_s|+|w_{ss}|+|w_z| + |w_{zz}|\big)(z,\xi)
 \le C(s,z) (1+|\xi|),&
\label{e.est1}
 \\
 &\big(|w_\xi| + |w_{z\xi}| + |w_{s\xi}|\big)(s,z)
 \le 
 C(s,z)
 ~\mbox{and}~
 |w_{\xi\xi}|
 \le
 \big(C\1_{[\xi_0,\xi_1]}\big)(s,z),&
\label{e.est2}
 \ee
where $C$ is an appropriate continuous function in $\R_+^2$,
depending on the Merton value function and its derivatives.}
\qed
\end{Remark}

\subsection{Remainder Estimate}
\label{s.remainder}

\label{ss.rem}

In this subsection, we estimate the remainder term 
in Remark \ref{r.formal}. So, let $\Psi^\ep$ be as in Remark \ref{r.formal} with $\varpi$ satisfying the same estimates \reff{e.est1}-\reff{e.est2} as $w$. We have seen in \reff{e.rem} that
\begin{eqnarray*}
\Jc(\Psi^\ep)(s,x,y)
&:=& 
\big(\beta\Psi^\ep-\Lc \Psi^\ep-\Ut(\Psi^\ep_x)\big)(s,x,y)\\
&=&
\ep^2 \left[-\frac12 v_{zz}(s,z) \xi^2 
            +\frac12\alpha^2(s,z) \varpi_{\xi \xi}(s,z,\xi) 
            -\Ac \phi (s,z)+\Rc^\ep(s,z,\xi)
      \right],
 \end{eqnarray*}
where $\alpha$, $\Ac$ are defined in \reff{alpha}-\reff{Ac}, and $\Rc^\ep$ is the remainder. By a direct (tedious) calculation, 
the remainder term can be obtained explicitly. 
In view of our previous bounds \reff{e.est1}-\reff{e.est2} 
on the derivatives of $w$, we obtain the estimate,
 \begin{eqnarray*}
 \big|\Rc^\ep(s,z,\xi)\big|
 &\le&
 \eps\Big(|\xi| |\mu-r| |\phi_z| 
        +\frac12\sigma^2(\eps\xi^2+2|\xi||\ybf|)|\phi_{zz}|
        +\sigma^2|\xi||\phi_{sz}|
      \Big)(s,z) 
  \\
  &&
  +\eps C(s,z)\big(1+\eps|\xi|+\eps^2|\xi|^2+\eps^3|\xi|^3\big),
 \\
 &&
 +\eps^{-2}
 \big|\tilde U(\psi_x^\eps)-\tilde U(v_z)
      -(\psi_x^\eps-v_z)\tilde U'(v_z)
 \big|
 \end{eqnarray*}
for some continuous function $C(s,z)$. 
Since $\tilde U$ is $C^1$ and convex, 
 \begin{eqnarray*}
 \big|\Rc^\ep(s,z,\xi)\big|
 &\le&
 \eps\Big(|\xi| |\mu-r| |\phi_z| 
        +\frac12\sigma^2(\eps\xi^2+2|\xi||\ybf|)|\phi_{zz}|
        +\sigma^2|\xi||\phi_{sz}|
      \Big)(s,z) \\
  &&
  +\eps C(s,z)\big(1+\eps|\xi|+\eps^2|\xi|^2+\eps^3|\xi|^3\big),
 \\
 &&
 +(|\phi_z|+\eps^2 |\phi_z|+\eps\ybf_z|\varpi_\xi|)
  \big|\tilde U'(v_z)+\eps^2|\phi_z|+\eps^4 |\varpi_z|
                    +\eps^3\ybf_z|\varpi_\xi|)
      -\tilde U'(v_z)
 \big|
 \e*
Suppose that $\varpi$ satisfies the same estimates 
 \reff{e.est1}-\reff{e.est2} as $w$.
 Then,
 \b*
 \big|\Rc^\ep(s,z,\xi)\big|
 &\le&
 \eps\Big(|\xi| |\mu-r| |\phi_z| 
        +\frac12\sigma^2(\eps\xi^2+2|\xi||\ybf|)|\phi_{zz}|
        +\sigma^2|\xi||\phi_{sz}|
      \Big)(s,z)
  \\
  &&
  +\eps C(s,z)\big(1+\eps|\xi|+\eps^2|\xi|^2+\eps^3|\xi|^3\big),
 \\
 &&
 +\eps^2\big(|\phi_z|+\eps C(s,z)(1+\eps|\xi|)\big)^2
  \tilde U''\big(v_z+\eps^2|\phi_z|+\eps^3C(s,z)(1+\eps|\xi|)\big).
 \end{eqnarray*}
\vspace{1pt}

\section{Assumptions}
\label{s.assume}

The main objective of this paper is to characterize the limit of
the following sequence,

$$
 \bar{u}^\ep(s,x,y)
 :=
 \frac{v(s,z)-\ve(s,x,y)}
      {\eps^2},
 ~~s\ge 0,~(x,y)\in K_\eps.
$$

Our proof follows the general methodology developed by 
Barles \& Perthame in the context of viscosity solutions. 
Hence, we first define relaxed semi-limits by, 

 \be
\nonumber
 u^*(\zeta)
 := \limsup_{(\ep,\zeta^\prime) \to (0,\zeta)} \   
     \bar{u}^\ep(\zeta^\prime),
 &&
 u_*(\zeta)
 := \liminf_{(\ep,\zeta^\prime) \to (0,\zeta)} \
     \bar{u}^\ep(\zeta^\prime)
 \ee
Then, we show  under appropriate conditions that they 
are viscosity sub-solution and super-solution, respectively, 
of the second corrector equation 
\reff{e.au}.

We shall now formulate some conditions which guarantee that
\begin{itemize}
\item[i.] the relaxed semi-limits are finite, 
\item[ii.] the second corrector equation \reff{e.au} 
verifies comparison for viscosity solutions.
\end{itemize}
We may then conclude that $u^*\le u_*$. 
Since the opposite inequality is obvious, 
this shows that $u=u^*=u_*$ is 
the unique solution of the second corrector equation \reff{e.au}.

In this short subsection, for the convenience of 
the reader, we collect all the assumptions 
needed for the convergence proof, including 
the ones that were already used. 

We first focus on the finiteness of the relaxed
 semi-limits $u_*$ and $u^*$. 
A local lower bound is easy to obtain in view of the obvious inequality 
$v^\eps(s,x,y)\le v(s,x+y)$ which implies that $\bar{u}^\ep \ge 0$.
Our first assumption complements this with a local upper bound.

\begin{Assumption}[Uniform Local Bound]
\label{a.bound}
The family of functions $\bar{u}^\ep$ is locally uniformly
bounded from above.
\end{Assumption}

The above assumption states that 
for any $(s_0,x_0,y_0) \in \R_+\x\R^2$ with $x_0+y_0>0$, 
there exist $r_0=r_0(s_0,x_0,y_0)>0$ and $\ep_0=\ep_0(s_0,x_0,y_0)>0$ so that
 \begin{equation} 
 \label{e.bound}
 b(s_0,x_0,y_0):= 
 \sup\{\ \ue(s,x,y)\ :\ (s,x,y) \in B_{r_0}(s_0,x_0,y_0),\
 \ep \in (0,\ep_0]\ \} <\infty,
 \end{equation}
where $B_{r_0}(s_0,x_0,y_0)$ denotes the open ball with radius $r_0$, centered at $(s_0,x_0,y_0)$.

This assumption is verified in Section \ref{s.verify} under 
some conditions on $v$ and its derivatives by 
constructing an appropriate 
sub-solution to the dynamic programming equation \reff{e.dpp}. 
However, the sub-solution does not need to have the exact $\ep^2$ behavior 
as needed in other approaches to this problem starting from \cite{ss,js2004}. 
Indeed, in these earlier approaches, both the sub and the super-solution
must be sharp enough to have the exact limiting behavior in 
the leading $\ep^2$ term. For the above estimate, however, 
this term needs to be only locally bounded.

The next assumption is a  regularity condition on the Merton problem.

\begin{Assumption}[Smoothness]
\label{a.smooth} The Merton value function $v$ and the 
Merton optimal investment strategy $\ybf$ are twice continuously 
differentiable in the open domain $(0,\infty)^2$ and $v_z(s,z)>0 $ for all $s,z>0$. 
Moreover, there exist $c_1\ge c_0>0$ such that 
\be
\label{a.theta}
c_0 z\le [\ybf(1-\ybf_z)-s\ybf](s,z) \le c_1 z
&\mbox{for all}&
s,z \in \R_+.
\ee
\end{Assumption}

In particular, together with our condition standing assumption 
on the volatility function $\sigma$, the above assumption implies 
that the diffusion coefficient $\alpha(s,z)$ in the first corrector 
equation is  non-degenerate away from the origin. 
For later use we record that there exist two constants $0<\alpha_*\le \alpha^*$ so that
\begin{equation}
\label{e.alphabound}
0< \alpha_* \le
\frac{\alpha(s,z)}{z} \le \alpha^*,
\qquad
\forall \ s,z \in \R_+.
\end{equation}
We will not attempt to verify the above hypothesis. 
However, in the power utility case, the value function is always 
smooth and the condition \reff{a.theta} can be directly checked 
as the optimal investment policy $\ybf$ is explicitly available.
 
We next assume that the second corrector equation \reff{e.au}
 has comparison. Recall the function $u$ introduced in \reff{e.defu}, 
 let $b$ be as in \reff{e.bound}, and set
\begin{equation}
\label{e.B}
B(s,z):= b\big(s,z-\ybf(z), \ybf(z)\big),
\qquad
s,z \in \R_+.
\end{equation}

\begin{Assumption}[Comparison]
\label{a.compare}
For any upper-semicontinuous {\rm{(}}resp.~lower-semicontinuous{\rm{)}} viscosity sub-solution 
{\rm{(}}resp.~super-solution{\rm{)}}
$u_1$ {\rm{(}}resp. $u_2${\rm{)}} of \reff{e.au} in $(0,\infty)^2$ satisfying the growth condition
$|u_i| \le B$ on $(0,\infty)^2$, $i=1,2$, we have $u_1 \le u \le  u_2$ in $(0,\infty)^2$.
\end{Assumption}

In the above comparison, notice that the growth of the supersolution 
and the subsolution is controlled by the function $B$ which is defined 
in \reff{e.B} by means of the local bound function $b$. In particular,
$B$ controls the growth both at infinity and near the origin. This observation is further detailed in Remark \ref{rem:behaviorzero} below. 

We observe however that, as discussed earlier, the operator 
$\Ac$ is the infinitesimal generator of the optimal wealth process
 in the limiting Merton problem. In view of our Assumption \ref{a.smooth}, 
 we implicitly assume that this process does not reach the origin with probability one. 

We finally formulate a natural assumption which was verified in 
\cite{ss}, Remark 11.3, in the context of power utility functions. 
This assumption will be used for the proof of the sub-solution property. 
To state this assumption, we first introduce the {\em{no-transaction region}} 
defined by,   
 \begin{equation}
\label{e.notrade}
 \Nc^\ep
 :=
 \left\{(s,x,y) \in K_\ep\ : \Lambda^\ep_{0,1}\cdot D\ve(s,x,y)>0, \ 
                             \mbox{and} \
                             \Lambda^\ep_{1,0}\cdot D\ve(s,x,y)>0
\right\}.
 \end{equation}
By the dynamic programming equation \reff{e.dpp}, the value function $v^\eps$ is a viscosity solution of 
 \b*
 \beta\ve-\Lc \ve-\tilde U(\ve_x) =0
 &\mbox{on}&
 \Nc^\ep.
 \e*
 
\begin{Assumption}[No transaction region] \label{a.merton} 
The no-transaction region $\Nc^\ep$ contains the Merton line 
$\Mc:= \{ (s,z-\ybf(z), \ybf(z))\ :\ s,z \in \R_+\ \}$.
\end{Assumption} 

\begin{Remark}{\rm
In our accompanying paper \cite{pst}, the expansion result in 
the $d-$dimensional context is proved without Assumption \ref{a.merton}. 
However, this induces an important additional technical effort. 
Therefore, for the sake of simplicity, we refrained from including 
this improvement in the present one-dimensional paper.
}
\end{Remark}

\section{Convergence in one dimension}
\label{s.convergence}
For the convergence proof, we introduce the 
following ``corrected'' version of $\bar{u}^\ep$,
$$
 {u}^\ep(s,x,y)
 := \bar u^\eps(s,x,y)-\ep^2 w(s,z,\xi),
 ~~s\ge 0,~(x,y)\in K_\eps.
$$
 Notice that both families
 $\bar{u}^\ep$ and $u^\ep$ have the same
 relaxed semi-limits $u^*$ and $u_*$.

\begin{Theorem}
\label{t.main}
Under Assumptions \ref{a.bound}, \ref{a.smooth}, 
\ref{a.compare}, and \ref{a.merton}
the sequence $\{u^\ep\}_{\eps>0}$ converges locally 
uniformly to the function $u$ defined in \reff{e.defu}.
\end{Theorem}

\proof
In the next subsections, we will show that, 
 the semi-limits
$u_*$ and $u^*$ are viscosity super-solution and sub-solution, 
respectively, of \reff{e.au}.
Then, by the comparison Assumption \ref{a.compare}, we 
conclude
that $u^* \le u \le  u_*$.  Since the opposite inequality
is obvious, this implies that $u^*=u_*=u$.
The local uniform convergence follows immediately
from this and the definitions.
\qed

\subsection{First properties}

In this subsection, we only use the assumptions on the smoothness 
of the limiting Merton problem and the local boundedness of 
$\{u^\eps\}_{\eps}$. We first recall that
 \b*
 \overline{\lambda}
 &:=&
 \lambda^{0,1}\vee\lambda^{1,0}.
 \e*

\begin{Lemma}\label{l.locbound}
{\rm (i)} For all $\eps,s>0$, $(x,y)\in K_\eps$, 
$u^\eps(s,x,y)\ge -\eps \overline{\lambda}v_z(s,z)|y-\ybf(s,z)|$. 
In particular, $u^*\ge 0$.
\\
{\rm (ii)} If in addition Assumption \ref{a.bound} holds, then 
 \b*
 0 \le u_*(s,x,y)\le u^*(s,x,y)<\infty
 &\mbox{for all}& 
 s,x,y>0.
 \e*
\end{Lemma}

\no {\em{Proof.}} Since (ii) is a direct consequence, 
we focus on (i). From the obvious inequality 
$\ve(s,x,y) \le v(s,x+y)$, it follows that $\ue(s,x,y) \ge -\ep^2 w(s,z,\xi)$, 
so that the required result follows from the bound 
\reff{bound-wxi} on $w_\xi$ together with $w(\cdot,0)=0$. 
\qed

\vspace{5mm}

We next show that the relaxed semi-limits  $u^*$ and $u_*$ 
depend on the pair $(x,y)$ only through the aggregate variable 
$z=x+y$. 

\begin{Lemma}
\label{l.reduction}
Let Assumptions \ref{a.bound} and \ref{a.smooth} hold true. 
Then, $u^*$ and $u_*$ are functions of $(s,z)$ only. 
Moreover, for all $s,z\ge 0$,
 $$
 u_*(s,z)
 =\!\!\liminf_{(\eps,s',z')\to(0,s,z)} 
  u^\eps\big(s',z'-\ybf(z'),\ybf(z')\big),
$$
and
$$
 u^*(s,z)
 =\!\!\limsup_{(\eps,s',z')\to(0,s,z)} 
  u^\eps\big(s',z'-\ybf(z'),\ybf(z')\big).
 $$
\end{Lemma}

\noindent {\em{Proof.}} This result is a consequence 
of the gradient constraints in the 
dynamic programming equation \reff{e.dpp},
$$
 \Lambda^\eps_{1,0}\cdot (v^\eps_{x},v^\eps_{y})\ge 0
 ~~\mbox{and}~~
 \Lambda^\eps_{0,1}\cdot  (v^\eps_{x},v^\eps_{y}) \ge 0
 \ \ \mbox{in the viscosity sense.}
$$
{\bf 1.} We change variables and use  the above inequalities to obtain
 \be\label{indep-z1}
 \big(1+ \lambda^{1,0}\eps^3(1-\ybf_z)\big)\hat v^\eps_\xi
 \ge 
 - \lambda^{1,0}\eps^4\hat v^\eps_z,
 &&
 \big(1+\lambda^{0,1}\eps^3 \ybf_z \big)\hat v^\eps_\xi
 \le 
 \lambda^{0,1}\eps^4\hat v^\eps_z,
 \ee
in the viscosity sense. Since $v^\eps$ is concave in $(x,y)$, 
the partial gradients $v^\eps_x$ and $v^\eps_y$ exist almost everywhere. 
By the smoothness of the Merton optimal investment strategy 
$\ybf$, this implies that the partial gradient $\hat v^\eps_z$ also exists almost everywhere. 
Then, by the definition of $u^\eps$, 
we conclude that the partial gradients 
$\hat u^\eps_z$ and $\hat u^\eps_\xi$ exist almost everywhere. 
In view of Condition \reff{a.theta} in Assumption \ref{a.smooth}, 
we conclude from \reff{indep-z1} and the fact that $\hat v^\eps_z\ge 0$ that
\be
\label{hatvxi}
\left| \hat v^\ep_\xi\right|
\le \overline{\lambda} \eps^4\hat v^\eps_z.
\ee
We now claim that
 \begin{eqnarray}
\nonumber
 \hat v^\eps_z(s,z,\xi)
& \le&
 \gamma^\eps(s,x,y)
 \\
 &&
 :=v_z(s,z-\ep) +
 \eps \big(u^\eps(s,x-\eps,y)+u^\eps(s,x,y-\eps)\big)
 \label{claim-hatvz}
 \\
 && +\eps^3\overline{\lambda}v_z(s,z) 
  \Big(1+|\ybf_z(s,z)|+|\xi|+\frac{|\ybf(s,z)-\ybf(s,z-\eps)|}
                    {\eps}
  \Big).
 \nonumber
 \end{eqnarray}
We postpone the justification of this claim
to the next step and continue with the proof. 
Then, it follows from \reff{hatvxi}, \reff{claim-hatvz} 
together with Assumption \ref{a.smooth} and 
\reff{bound-wxi},
 \be
 \big|\hat u^\eps_\xi(s,z,\xi)\big|
 &\le&
 \eps^2 \bar{\lambda}  \left( v_z(s,z)
 +\hat v^\eps_z(s,z,\xi)\right)
 \nonumber\\
 &\le&
 \eps^2 \bar{\lambda}  \left( v_z(s,z)
 +  \gamma^\eps(s,z,\xi)\right).
 \label{hatuxi}
 \ee
Hence,
$$
 (e_1-e_0)\cdot (u^\ep_x,u^\ep_y)=
 -\frac{1}{\ep} \hat u^\eps_\xi\le
 \eps \bar{\lambda}  \left( v_z(s,z)
 +  \gamma^\eps(s,z,\xi)\right).
$$
By the  local boundedness of $\{u^\eps\}_{\eps}$,
for any $(s,x,y)$, there is an open neighborhood 
of $(s,x,y)$ and a constant $K$, 
both independent of $\eps$, such that the maps
 \b*
 t\mapsto u^\eps(s,x-t,y+t)+\eps Kt
 &\mbox{and}&
 t\mapsto -u^\eps(s,x-t,y+t)+\eps Kt
 \e*
are nondecreasing for all $\eps>0$. Then, 
it follows from the definition of the relaxed semi-limits 
that $\hat u^*$ and $\hat u_*$ are independent of the $\xi$-variable.

\noindent
{\bf 2.} We now prove \reff{claim-hatvz}. 
For $\eps>0$ and $(x,y), (x-\eps,y), (x,y-\eps)\in K_\eps$, 
we denote as usual $z=x+y$ and $\xi=(y-\ybf(s,z))/\eps$. 
By the concavity of $v^\eps$ in the pair $(x,y)$ and 
the concavity of the Merton function $v$ in $z$ that:
 \b*
 v^\eps_x(s,x,y)
 &\le&
 \frac{1}{\eps}\big(v^\eps(s,x,y)-v^\eps(s,x-\eps,y)\big)
 \\
 &\le&
 \frac{1}{\eps}\big(v(s,z)-v(s,z-\eps)\big)
 +\frac{1}{\eps}\big(v(s,z-\eps)-v^\eps(s,x-\eps,y)\big)
 \\
 &\le&
 v_z(s,z-\eps)
 +\frac{1}{\eps}\big(v(s,z-\eps)-v^\eps(s,x-\eps,y)\big).
 \e*
By the definition of $u^\eps$,
$$
 v^\eps_x(s,x,y)
 \le
 v_z(s,z-\eps)
 +\eps\big(u^\eps(s,x-\eps,y)+\eps^2w(s,z-\eps,\xi_\eps)\big)
$$
where $\xi_\eps:=(y-\ybf(s,z-\eps))/\eps=
\xi+(\ybf(s,z)-\ybf(s,z-\eps))/\eps$. 
We use the bound \reff{boundwwxi} on $w$, to arrive at,
 $$
 v^\eps_x(s,x,y)
 \le v_z(s,z-\ep) +
\eps u^\eps(s,x-\eps,y)
 +\eps^3\overline{\lambda}v_z(s,z) 
  \Big(1+|\xi|+\frac{|\ybf(s,z)-\ybf(s,z-\eps)|}
                    {\eps}
  \Big).
$$
By exactly the same argument, we also conclude that
 $$
 v^\eps_y(s,x,y)
 \le v_z(s,z-\ep) +
\eps u^\eps(s,x,y-\eps)
 +\eps^3\overline{\lambda}v_z(s,z) 
  \Big(1+|\xi|+\frac{|-\eps+\ybf(s,z)-\ybf(s,z-\eps)|}
                    {\eps}
  \Big).
 $$
Then, using the bounds on $\ybf_z$ from Assumption \ref{a.smooth},
 \b*
 \hat v^\eps_z(s,z,\xi)
 &=&
 \partial_zv^\eps
 \big(s,z-\eps\xi-\ybf(s,z),\eps\xi+\ybf(s,z)\big)
 \\
 &=&
 \big(1-\ybf_z(s,z)\big)v^\eps_x(s,x,y)
 +\ybf_z(s,z)v^\eps_y(s,x,y)
 \\
 &\le& v_z(s,z-\ep) +
 \eps \big(u^\eps(s,x-\eps,y)+u^\eps(s,x,y-\eps)\big)
 \\
 &&
 +\eps^3\overline{\lambda}v_z(s,z) 
  \Big(1+|\ybf_z(s,z)|+|\xi|
       +\frac{|\ybf(s,z)-\ybf(s,z-\eps)|}
             {\eps}
  \Big).
 \e*
{\bf 3.} The final statement in the lemma follows 
from \reff{hatuxi}, the expression of $\gamma^\eps$ in \reff{claim-hatvz}, 
and Assumption \ref{a.bound}.
\qed

\subsection{Viscosity sub-solution property}

In this section, we prove

\begin{Proposition}\label{prop-supersol}
Under Assumptions \ref{a.bound} and \ref{a.smooth},
 the function $u^*$ is a viscosity sub-solution of 
the second corrector equation \reff{e.au}.
\end{Proposition}

\noindent {\em {Proof.}} 
Let $(s_0,z_0,\varphi)\in(0,\infty)^2\times C^2(\R_+^2)$ be such that
 \be
 \label{strictmin}
 0=(u^*-\varphi)(s_0,z_0)
 >
 (u^*-\varphi)(s,z)
 &\mbox{for all}&
 s,z\ge 0,~~(s,z)\neq (s_0,z_0).
 \ee
Our objective in the following steps is to prove that
 \bee
 \label{goal-supersol}
 \Ac\varphi(s_0,z_0) - a(s_0,z_0) \le 0.
 \eee

\no {{\bf{1.}} 
By the definition of $u^*$ and Lemma \ref{l.reduction}, 
there exists a sequence $(s^\eps,\ze)$ so that
 \b*
 (s^\eps,\ze)\to (s_0,z_0)
 ~~\mbox{and}~~ 
 \hat{u}^\ep(s^\eps,\ze,0) 
 \to
 u^*(s_0,z_0),
 &\mbox{as}~~\eps\downarrow  0,&
 \e*
where we used the notation \reff{e.hat}. Then, it is clear that
 \bee
 \label{ceps}
 \ell^\eps_* :=\hat{u}^\ep(s^\eps,\ze,0)-\varphi(s^\eps,\ze)
 \to 0
 \eee
and
\b*
(x^\ep,y^\ep)=\big(z^\ep-\ybf(s^\eps,z^\ep),\ybf(s^\eps,z^\ep)\big)
&\longrightarrow& 
(x_0,y_0):=\big(z_0-\ybf(s_0,z_0),\ybf(s_0,z_0)\big).
\e*
Since $(u^\eps)$ is locally bounded 
from above (Assumption \ref{a.bound}), 
there are
$r_0:=r_0(s_0,x_0,y_0)>0$ and 
$\ep_0:=\ep_0(s_0,x_0,y_0)>0$ so that
 \be
 \label{b*}
 b_*:=\sup \{ \ue(s,x,y)\ :\
 (s,x,y)\in B_0, \ep \in (0,\ep_0]\}<\infty,
 \quad \mbox{where}\quad
 B_0:=B_{r_0}(s_0,x_0,y_0)
 \ee 
is the open ball centered at $(s_0,x_0,y_0)$ with radius $r_0$. 
We may choose $r_0\le z_0/2$ so that $B_0$ 
does not intersect the line $z=0$. 
For $\eps, \delta \in(0,1]$, set 
 \b*
 \hat{\psi}^{\eps,\delta}(s,z,\xi)
 &:=& 
 v(s,z)-\ep^2 \ell^\eps_*
 -\eps^2\varphi(s,z)
 -\eps^4 (1+ \delta)  w(s,z,\xi)
 - \eps^2 \hat \phi^\eps(s,z,\xi),
 \e*
where, following our standard notation \reff{e.hat}, 
$\hat\phi^\eps$ is determined from the function,
 \b*
 \phi^\ep(s,x,y)
 &:=& 
 C\left[ (s-s^\eps)^4+(x+y-\ze)^4 + (y-\ybf(s,x+y))^4\right],
 \e*
and $C>0$ is a large constant that is chosen so that for all sufficiently
small $\ep>0$,
 \be\label{e.varphi}
 \phi^\ep\ge 1+ b_*-\varphi,
 ~~\mbox{on}~~ 
 B_0\setminus B_1
 &\mbox{with}&
 B_1:=B_{r_0/2}(s_0,x_0,y_0).
 \ee
The constant $C$ chosen above may depend on 
many things including the test function 
$\varphi$, $s_0,z_0,\delta$, but not on $\epsilon$. 
The convergence of $(s^\eps,\ze)$ to $(s_0,z_0)$ 
determines how small $\ep$ should be for \reff{e.varphi} to hold.

\no {{\bf{2.}} We first show that, for all sufficiently small
$\ep>0$, $\delta>0$, the difference 
$(v^\eps-\psi^{\eps,\delta})$, or equivalently,
\begin{eqnarray*}
 I^{\ep,\delta}(s,x,y)
& :=&
 \frac{v^\eps(s,x,y)-\psi^{\eps,\delta}(s,x,y)}
      {\ep^2}\\
& =&
 -u^\eps(s,x,y)+\varphi(s,z)+\ell^\eps_*+ \phi^\ep(s,x,y)
            +\eps^2\delta w(s,z,\xi),
 \end{eqnarray*}
has a local minimizer in $B_0$. Indeed, by the definition of 
$\ue$, $\psi^{\eps,\delta}$ and $\ell^\eps_*$, 
\reff{e.varphi}, \reff{b*}, and the fact that $w\ge 0$ that,
for any $(s,x,y)\in\partial B_0$,
 \b*
 I^{\ep,\delta}(s,x,y)
 &\ge&
 -u^\eps(s,x,y)+\ell^\eps_*+1+b_*+\eps^2\delta w(s,z,\xi)
 \;\ge\;
 1+\ell^\eps_*
 \;>\;0,
 \e* 
for sufficiently small $\eps$ 
in view of  \reff{ceps}. 
Since $I^{\ep,\delta}(s^\eps,x^\ep,y^\eps)=0$, 
we conclude that $I^{\ep, \delta}$ has a local 
minimizer $(\tilde{s}^\eps,\tilde{x}^\eps,\tilde{y}^\eps)$ in 
$B_0$ with $\tilde z^\ep:= \tilde{x}^\eps+\tilde{y}^\eps$,
$\tilde \xi^\ep:= (\tilde{y}^\eps-\ybf(\tilde s^\eps,\tilde z^\ep))/\ep$ 
satisfying,
 $$
 \min_{(s,z,\xi)\in B_1}(\hat{v}^\eps-\hat{\psi}^{\eps,\delta})
 =
 (\hat{v}^\eps-\hat{\psi}^{\eps,\delta})(\tilde z_\eps,\tilde\xi_\eps) \le 0,
 \quad
 |\tilde s^\ep -s_0|+|\tilde z^\ep -z_0| < r_0,\ \
 |\xi_\ep| < r_1/\ep,
$$
for some constant $r_1$. Since $v^\eps$  is a viscosity 
super-solution of the dynamic programming equation 
\reff{e.dpp}, we conclude that
 \be \label{supersol-PDE}
  \left(\beta \ve
        - \Lc \psi^{\eps,\delta}-\tilde U\big(\psi^{\eps,\delta}_x\big)
  \right)(\tilde s^\eps,\tilde x^\eps,\tilde y^\eps)      
  &\ge& 
  0,
 \ee
and 
 \b*
 \Lambda^\eps_{1,0}
 \cdot 
 \big(\psi^{\eps,\delta}_x,\psi^{\eps,\delta}_y\big)
 (\tilde s^\eps,\tilde x^\eps,\tilde y^\eps)
 &=&
 \left(\psi^{\eps,\delta}_x
       -(1- \lambda^{1,0}\ep^3)\psi^{\eps,\delta}_y\right)
                               (\tilde s^\eps,\tilde x^\eps,\tilde y^\eps) 
 \;\ge\; 0,
 \\
 \Lambda^\eps_{0,1}\cdot 
 \big(\psi^{\eps,\delta}_x,\psi^{\eps,\delta}_y\big)
 (\tilde s^\eps,\tilde x^\eps,\tilde y^\eps)
 &=&
 \left(\psi^{\eps,\delta}_y
       -(1-\lambda^{0,1}\ep^3)\psi^{\eps,\delta}_x\right)
                              (\tilde s^\eps,\tilde x^\eps,\tilde y^\eps) 
 \;\ge\; 0.
 \e*
By a direct calculation using the boundedness of 
$(\tilde s^\eps,\tilde z^\eps,\eps\tilde\xi^\eps)$, 
we rewrite the last gradient inequalities as follows,
 \be
 -4\eps^2(\eps\tilde\xi^\eps)^3
 +\eps^3 v_z(\tilde s^\eps,\tilde z^\eps)
         \big[\lambda^{1,0}
              -(1+\delta)
               \overline{w}_\rho(\tilde s^\eps,\tilde z^\eps,\tilde\rho^\eps)
         \big]
 +\circ(\eps^3)
 &\ge& 0,
 \label{supersol-nabla10}\\
 4\eps^2(\eps\tilde\xi^\eps)^3
 +\eps^3 v_z(\tilde s^\eps,\tilde z^\eps)
         \big[\lambda^{0,1}
              +(1+\delta)
               \overline{w}_\rho(\tilde s^\eps,\tilde z^\eps,\tilde\rho^\eps)
         \big]
 +\circ(\eps^3)
 &\ge& 0,
 \label{supersol-nabla01}
 \ee
where $\tilde\rho^\eps:=\tilde\xi^\eps
                        /\eta(\tilde{s}^\eps,\tilde{z}^\eps)$.

\no {{\bf{3.}} Let $\rho_0(s,z)$ be as in \reff{e.xi}. In this step, we show that 
 \be\label{goal-step3}
 |\tilde\rho^\eps|<\rho_0(\tilde s^\eps,\tilde z^\eps)
 &\mbox{for all sufficiently small}&
 \eps\in(0,1].
 \ee 
Indeed, assume that 
$\tilde\rho^{\eps_n} \le -\rho_0(\tilde s^{\eps_n},\tilde z^{\eps_n})
=\rho_1(\tilde s^{\eps_n},\tilde z^{\eps_n})$ for some sequence 
$\eps_n \in(0,1]$ with $\eps_n\to 0$. 
Then, 
$\overline{w}_\rho(\tilde s^{\eps_n},\tilde z^{\eps_n},\tilde \rho^{\eps_n})
=-\lambda^{0,1}$, and it follows from inequality 
\reff{supersol-nabla01}, together with the fact 
$\tilde\rho^{\eps_n}\le\rho_1(\tilde s^{\eps_n},\tilde z^{\eps_n})\le 0$, 
that
 $$
 0 
 \le
 4\eps_n^2({\eps_n}\tilde\xi^{\eps_n})^3
 -\eps_n^3 v_z(\tilde s^{\eps_n},\tilde z^{\eps_n})\delta\lambda^{0,1}
 +\circ(\eps_n^3)
 \le
 -{\eps_n}^3 v_z(\tilde s^{\eps_n},\tilde z^{\eps_n})\delta\lambda^{0,1}
 +\circ({\eps_n}^3).
 $$  
Since $\delta>0$, this can not happen for large $n$. 
Similarly, if 
$\tilde\rho^{\eps_n} \ge \rho_0(\tilde s^{\eps_n},\tilde z^{\eps_n})$ 
for some sequence $\eps_n\to 0$, we have 
$\overline{w}_\rho(\tilde s^{\eps_n},\tilde z^{\eps_n},\tilde \rho^{\eps_n})
=\lambda^{1,0}$, and it follows from inequality \reff{supersol-nabla10}, 
together with the fact that 
$\tilde\rho^{\eps_n}\ge\rho_0(\tilde s^{\eps_n},\tilde z^{\eps_n})\ge 0$, that
 $$
 0 
 \le
 -4\eps_n^2({\eps_n}\tilde\xi^{\eps_n})^3
 +\eps_n^3 v_z(\tilde s^{\eps_n},\tilde z^{\eps_n})(-\delta\lambda^{1,0})
 +\circ({\eps_n}^3)
 \le
 -\eps_n^3 v_z(\tilde s^{\eps_n},\tilde z^{\eps_n})\delta\lambda^{1,0}
 +\circ(\eps_n^3),
 $$ 
which leads again to a contradiction for large $n$, 
completing the proof of \reff{goal-step3}. 

\no {{\bf{4.}} Since $(\tilde s^\eps,
\tilde z^\ep)$ is bounded and 
$(s,z)\mapsto\rho_0(s,z)$ is continuous, 
we conclude from \reff{goal-step3} that the sequence 
$(\tilde\xi^\eps)_\eps$ is bounded. 
Hence, there exists a sequence $\eps_n\to 0$ so that
 \b*
 (s_n,z_n,\xi_n)
 :=(\tilde s^{\eps_n},\tilde z^{\eps_n},\tilde\xi^{\eps_n})
 &\longrightarrow& 
 (\hat s,\hat z,\hat\xi)=(s_0,z_0,\hat\xi)
 \e*
for some $\hat\xi\in\R$. The fact that the limit of $(s_n,z_n)$ is equal to
$(s_0,z_0)$ follows from standard arguments using the 
strict minimum property of $(s_0,z_0)$ in \reff{strictmin}. 
We now take the limit in \reff{supersol-PDE} 
along the sequence $\eps_n$. Since the function $\psi^{\ep,\delta}$ 
has the form as in Remark \ref{r.formal}, we do not repeat the 
computations given in Section \ref{s.formal} and, given 
the remainder estimate of section \ref{s.remainder}, we directly conclude that 
 \begin{eqnarray}
 \nonumber
 0
 &\le&
 \lim_{\ep_n \to 0}
 \ep_n^{-2}
 \Big(\beta v^{\eps_n} 
 - \Lc \psi^{\eps_n,\delta} 
 - \Ut\big(\psi^{\eps_n,\delta}_{x}
 \big)
 \Big)
 (s_n,z_n,\xi_n)
 \\
 &=&
 \frac12(\eta\sigma^2)(s_0,z_0)\hat\xi^2  
 +\frac12(1+\delta) \alpha^2(s_0,z_0) w_{\xi\xi}(s_0,z_0,\hat\xi) 
 - \Ac\varphi(s_0,z_0)
 \label{e.limit}
 \end{eqnarray}
In the above, we also used the fact that all derivatives of $\phi^\eps$ 
vanish at the origin as $\eps$ tends to zero.

\no {{\bf{5.}} In  Step 3, we have proved 
that $|\rho_\ep| \le \rho_0(z_\ep)$.
Hence, $|\hat\xi|\le (\eta\rho_0)(s_0,z_0)$. 
Since $w=\eta v_z \bar{w}$, $a =\eta v_z \bar{a}$,
the first corrector equation \reff{e.cw} implies that
 $$
 a(s_0,z_0)
 = 
 \frac12(\sigma^2\eta)(s_0,z_0)\hat\xi ^2
 +\frac12\alpha^2(s_0,z_0) w_{\xi\xi}(s_0,z_0,\hat\xi).
 $$ 
We use the above identity in \reff{e.limit}. The result is
 \b*
 \Ac\varphi(s_0,z_0)
 &\le& 
 \frac12(\sigma^2\eta)(s_0,z_0)\hat\xi ^2
 +\frac12(1+\delta)\alpha^2(s_0,z_0) w_{\xi\xi}(s_0,z_0,\hat\xi)
 \\
 &=& 
 a(s_0,z_0) 
 + \frac12\delta\alpha^2(s_0,z_0) w_{\xi\xi}(s_0,z_0,\hat\xi).
 \e*
Finally, we let $\delta$ go to zero. However, 
$\hat \xi=\hat\xi^\delta$ depends on $\delta$ 
and care must be taken. But since 
$|\xi_n| \le (\eta\rho_0)(s_n,z_n)$, it follows that 
$\hat\xi^\delta$ is uniformly bounded in $\delta$. 
Hence the second term in the above equation goes to 
zero with $\delta$, and we obtain the desired inequality 
\reff{goal-supersol}.
\qed

\subsection{Viscosity super-solution property}

In this section, we prove

\begin{Proposition}\label{prop-subsol}
Let Assumptions \ref{a.bound}, \ref{a.smooth}, and 
\ref{a.merton} hold true. 
Then, the function $u_*$ is a viscosity super-solution of the 
second corrector equation \reff{e.au}.
\end{Proposition}
 
As remarked earlier, the above result holds true
without the Assumption \ref{a.merton}
as proved in our forthcoming paper \cite{pst}. 
 However, in this paper we utilize it provide a somehow
 shorter proof.
We first need the following consequence of Assumption 
\ref{a.merton} and the convexity of $\ve$. 
Similar arguments are also used in \cite{ss}.

\begin{Lemma}\label{l.notrade}
Assume the hypothesis of Proposition \ref{prop-subsol}.
Let $(x,y)$ be an 
arbitrary element of $K_\ep$. Then,
\\
{\rm (i)} for $y \ge \ybf (s,z)$ {\rm{(}}or equivalently,
$\xi \ge 0${\rm{)}}, we have 
$\Lambda^\eps_{0,1}\cdot (\ve_x (s,x,y),\ve_y (s,x,y)) >0$,
\\
{\rm (ii)} for $y \le \ybf (s,z)$ {\rm{(}}or equivalently,
$\xi \le 0${\rm{)}}, we have
$\Lambda^\eps_{1,0}\cdot (\ve_x (s,x,y),\ve_y (s,x,y)) >0$.
\end{Lemma}

\no
{\em{Proof.}} For $z \in \R_+$ set
$$
\ybf^\ep_+(s,z):=
 \sup \big\{ y:  (z-y,y) \in K_\ep,\ \ 
                 {\mbox{and}}\ \ 
                 \Lambda_{0,1}^\eps\cdot (\ve_x,\ve_y )(s,z-y,y)=0 
      \big\}.
$$
In view of the form of $K_\ep$, we have 
$y \ge -z/(\eps^3\lambda^{0,1})$ and 
by convention the above supremum is equal 
to this lower bound if the set is empty. 
By the concavity of $\ve$, we conclude that
$$
 \Lambda_{0,1}^\eps\cdot (\ve_x,\ve_y )(s,x,y) 
 \left\{\begin{array}{l}
        =0~\mbox{for all}~ y \le \ybf^\ep_+(s,z),
        \\ 
        >0~\mbox{for all}~ y > \ybf^\ep_+(s,z).
        \end{array}
 \right.
$$
Let $\Nc^\eps$ be as in \reff{e.notrade}.  
Therefore it  is included in the set 
$\{ (s,x,y)\ :\ y > \ybf_+^\ep(s,z)\}$. 
Since Assumption \ref{a.merton} states that 
the Merton line $\{(s,x,y):~y=\ybf(s,z)\}$ is 
included in $\Nc^\eps$, we conclude that
$\ybf(s,z) > \ybf^\ep_+(s,z)$.
This proves the statement (i). 
The other assertion is proved similarly.
\qed

\vspace{5mm}

\noindent {\em {Proof of Proposition \ref{prop-subsol}}.} 
Let $(s_0,z_0,\varphi)\in(0,\infty)^2\times C^2(\R_+)$ be such that
 \be\label{strictmax}
 0=(u_*-\varphi)(s_0,z_0)
 <
 (u_*-\varphi)(s,z)
 &\mbox{for all}&
 s,z\ge 0,~~(s,z)\neq (s_0,z_0).
 \ee
We proceed to prove that
 \bee
 \label{goal-subsol}
 \Ac\varphi(s_0,z_0)-a(s_0,z_0) \ge 0.
 \eee

\no {{\bf{1.}}
By the definition of $u_*$ and Lemma \ref{l.reduction},
 there exists a sequence $(s^\eps,\ze)$ so that
 \b*
 (s^\eps,\ze)\to (s_0,z_0)
 ~~\mbox{and}~~ 
 \hat{u}^\ep(s^\eps,\ze,0) 
 \to
 u_*(s_0,z_0),
 &\mbox{as}~~\eps\downarrow 0,&
 \e*
where we used the notation \reff{e.hat}. Then, it is clear that
$$
 \ell_\eps^* :=\hat{u}^\ep(s^\eps,\ze,0)-\varphi(s^\eps,\ze)
 \longrightarrow 0
$$
and
$$
 (x^\ep,y^\ep)=\big(z^\ep-\ybf(s^\eps,z^\ep),\ybf(s^\eps,z^\ep)\big)
 \longrightarrow 
 (x_0,y_0):=\big(z_0-\ybf(s_0,z_0),\ybf(s_0,z_0)\big).
$$
Since $u^\eps(s,x,y)\ge-\eps^2w(s,z,\xi)\ge -\eps C(s,z)|y-\ybf(s,z)|$, 
for some continuous function $C$, there are
$r_0:=r_0(s_0,x_0,y_0)>0$ and $\ep_0:=\ep_0(s_0,x_0,y_0)>0$ so that
 \be
\nonumber
 b^*
 :=
 \inf_{(s,x,y)\in B_0} \ue(s,x,y)
 >-\infty,
 &\mbox{where}&
 B_0:=B_{r_0}(s_0,x_0,y_0).
 \ee 
We also choose $r_0$ sufficiently small so that $B_0$ does 
not intersect the line $z=0$. For $\eps\in(0,1]$ and $\delta>0$, define
$$
 \hat\psi^{\eps,\delta}(s,z,\xi)
 := v(s,z)
 -\ep^2 \ell_\eps^* 
 -\eps^2\varphi(s,z)
 -\eps^4 (1-\delta) w(s,z,\xi)
 +\eps^2\hat\phi^\eps(s,z,\xi),
$$
where, following our notation convention \reff{e.hat}, 
the function $\hat\phi^\eps$ is obtained from
the function $\phi^\eps$ defined by,
$$
 \phi^\eps(s,x,y)
 :=
 C\big[(s-s^\eps)^4+(x+y-z^\eps)^4+(y-\ybf(s,x+y))^4\big]
  \big]
$$
and, similar to the proof of the super-solution property, 
$C>0$ is a constant chosen so that,
 \be\label{C-subsol}
 -b^*+\ell^*_\eps
 +\big(\varphi-\phi^\eps\big)(s,x,y)
 <0
 &\mbox{on}&
 \partial B_0.
 \ee 
\no {{\bf{2.}} Set
 \b*
 I^{\eps,\delta}(s,z,\xi)
 &:=&
 \eps^{-2}\big(v^\eps-\psi^{\eps,\delta}\big)(s,x,y)
 \\
 &=&
 -u^\eps(s,x,y)+\varphi(s,z)+\ell_\ep^*
 -\phi^\eps(s,x,y)-\eps^2 \delta w(s,z,\xi).
 \e*
Since $w(s,z,0)=0$, we have $I^{\eps,\delta}(s^\ep,z^\ep,0)=0$. 
On the other hand, it follows from \reff{C-subsol} that 
 \b*
 I^{\eps,\delta}(s,z,\xi)
 \le
 -b^*+\ell^*_\eps
 +\big(\varphi-\phi^\eps\big)(s,x,y)
 -\eps^2 \delta w(s,z,\xi)
 <0
 &\mbox{on}&
 \partial B_0.
 \e*
Then, the difference $\ve-\psi^{\ep,\delta}$ 
has an interior maximizer 
$(\tilde s^\eps,\tilde z^\eps,\tilde\xi^\eps)$ in $B_0$,
 \begin{equation}
\label{zepsxieps}
 \max_{B_0}
 \big(v^\eps-\psi^{\lambda,\eps}\big)
 =
 (v^\eps-\psi^{\lambda,\eps})(\tilde s^\eps,\tilde x_\eps,\tilde y^\eps),
~\mbox{and}~
|\tilde s^\eps-s_0|
 +|\tilde z_\eps-z_0|
 +|\eps\tilde\xi_\eps|
 \le 
 r_1,
 \end{equation}
for some constant $r_1$. By the sub-solution property of $v^\eps$, 
at $(\tilde s^\eps,\tilde\xe,\tilde\ye)$,
 \begin{equation}
\label{psiepssubsol}
 \min\big\{\beta \ve
             - \Lc \psi^{\eps,\delta}
             -\tilde U\big(\psi^{\eps,\delta}_x\big),
           \Lambda^\eps_{0,1}\cdot (\psi^{\eps,\delta}_x,\psi^{\eps,\delta}_y),
           \Lambda^\eps_{1,0}\cdot (\psi^{\eps,\delta}_x,\psi^{\eps,\delta}_y)
     \big\}
 \le
 0.
 \end{equation}

\no {{\bf{3.}} In this step, we show that for all sufficiently small $\ep >0$, 
 \be \label{DL>0}
 \Lambda^\eps_{0,1}\cdot(\psi^{\eps,\delta}_x,\psi^{\eps,\delta}_y)
 (\tilde s^\eps,\tilde\xe,\tilde\ye)
 >0
 &\mbox{and}&
 \Lambda^\eps_{1,0}\cdot (\psi^{\eps,\delta}_x,\psi^{\eps,\delta}_y)
 (\tilde s^\eps,\tilde\xe,\tilde\ye)
 > 0.
 \ee
By Lemma \ref{l.notrade}, it suffices to prove that
 \begin{equation}
 \label{DL>0-1}
 \begin{array}{rcl}
D^{0,1}:= \Lambda^\eps_{0,1}\cdot (\psi^{\eps,\delta}_x,\psi^{\eps,\delta}_y)
 (\tilde s^\eps,\tilde\xe,\tilde\ye)
 >0
 &\mbox{for}&\tilde\xi<0,
 \\
D^{1,0}:= \Lambda^\eps_{1,0}\cdot (\psi^{\eps,\delta}_x,\psi^{\eps,\delta}_y)
 (\tilde s^\eps,\tilde\xe,\tilde\ye)
 > 0
 &\mbox{for}&
 \tilde\xi>0.
 \end{array}
 \end{equation}
We directly compute that
 \b*
\psi^{\eps,\delta}_x
 &=&
 v_z
 -\eps^2\varphi_z
 -\eps^4(1-\delta)\big(w_z-\frac{\ybf_z}{\eps}w_\xi\big)
 +4\eps^2 C\big((z-z^\eps)^3-\ybf_z(y-\ybf)^3
           \big),
 \\
\psi^{\eps,\delta}_y
 &=&
 v_z
 -\eps^2\varphi_z
 -\eps^4(1-\delta)\big(w_z+\frac{1-\ybf_z}{\eps}w_\xi\big)
 +4\eps^2 C\big((z-z^\eps)^3+(1-\ybf_z)(y-\ybf)^3
           \big).
 \e*
Then, it follows from the estimates \reff{zepsxieps} that
 \b*
 D^{0,1}
 &=&
 \eps^3\big((1-\delta)w_\xi+\lambda^{0,1}v_z\big)
       (\tilde s^\eps,\tilde z^\eps,\tilde\xi^\ep)
 -4C\eps^2(\eps\tilde\xi^\eps)^3
 +\circ(\eps^3)
 \\
D^{1,0}
 &=&
 \eps^3\big(-(1-\delta)w_\xi+\lambda^{1,0}v_z\big)
       (\tilde s^\eps,\tilde z^\eps,\tilde\xi^\ep)
 +4C\eps^2(\eps\tilde\xi^\eps)^3
 +\circ(\eps^3).
 \e*
 Since $w$ solves  \reff{e.1corrector-dim1},
 $w_\xi+\lambda^{0,1}v_z\ge 0$ and 
$-w_\xi+\lambda^{1,0}v_z\ge 0$.  Then,
\b*
 D^{0,1}
 &\ge&
 -\eps^3\delta v_z(\tilde s^\eps,\tilde z^\eps)
 -4C\eps^2(\eps\tilde\xi^\eps)^3
 +\circ(\eps^3)
 \\
 &\ge&
 -\eps^3\delta v_z(\tilde s^\eps,\tilde z^\eps)
 +\circ(\eps^3)
 ~~\mbox{for}~~\tilde\xi\le 0,
 \e*
and
 \b*
D^{1,0}
 &\ge&
 \eps^3\delta v_z(\tilde s^\eps,\tilde z^\eps)
 +4C\eps^2(\eps\tilde\xi^\eps)^3
 +\circ(\eps^3).
 \\
 &\ge&
 \eps^3\delta v_z(\tilde s^\eps,\tilde z^\eps)
 +\circ(\eps^3)
 ~~\mbox{for}~~\tilde\xi\ge 0.
 \e*
Since $v_z>0$, \reff{DL>0-1} holds for all sufficiently small $\eps>0$.

\no {{\bf 4.}}  In this step, we prove that 
$\tilde\xi_\eps$ is bounded in $\ep \in(0,1]$. 
Indeed, in view of  \reff{psiepssubsol} and \reff{DL>0}, 
 \be
 0
 &\ge&
 \Big(\beta v^\eps-\Lc\psi^{\eps,\delta}
      -\tilde U\big(\psi^{\eps,\delta}_x\big)
 \Big)(\tilde s^\eps,\tilde\xe,\tilde\ye)
 \nonumber\\
 &=&
 \nonumber
 \ep^2\Big[  \frac{(-\sigma^2 v_{zz})(s^\eps,\tilde z_\ep)}
                  {2}
             |\xi_\ep|^2
 +\frac{1-\delta}{2} \alpha^2(\tilde s^\eps,\tilde z_\ep) 
                     w_{\xi \xi}(\tilde z_\ep,\tilde\xi_\ep)\\
&&\qquad \qquad
  - \Ac u (\tilde s^\eps,\tilde z_\ep)
+ \Rc^\ep(\tilde s^\eps,\tilde x_\ep,\tilde y_\ep) \Big],
 \label{subsol-almostfinal}
 \ee
where we used the fact that the function $\psi^{\ep,\delta}$
 is exactly as in the form assumed in Remark \ref{r.formal}. 
Then, by the remainder estimate of section \ref{s.remainder}, 
we deduce that,
 \be\label{remainderestimate-subsol}
 |\Rc^\ep(\tilde s^\eps,\tilde x_\ep,\tilde y_\ep)| 
 \le
 C(\tilde s^\eps, \tilde z_\ep) 
 \left[ \ep +\ep|\tilde\xi_\ep|+ \ep^2 |\tilde\xi_\ep|^2
 \right].
 \ee
In Section \ref{s.corrector}, the function $w$ is explicitly
constructed. Since $w$ is linear in $\xi$ for large values of $\xi$,
there is a continuous function $\hat C(s,z)$ so that
 \b*
 0 \le w_{\xi \xi}(s,z,\xi) \le \hat C(s,z),
 &\mbox{for all}&
 (s,z,\xi) \in \R_+^2 \times \R^1.
 \e*
Then, since $(\tilde s^\eps,\tilde z_\ep)$ is 
uniformly bounded in $\ep \in (0,1]$, there are constants $C, \tilde C >0$ so that,
 $$ 0  \ge
 \ep^2\tilde C\left[\tilde\xi_\ep^2 
                       -C\left(1+\ep|\tilde\xi_\ep|+\ep^2|\tilde\xi_\ep|^2
                         \right)
                \right].
 $$
Hence $(\tilde\xi_\ep)_\ep$ is also uniformly 
bounded in $\ep \in (0,1]$ by a constant depending 
only on the test functions.

\no {{\bf{5.}}  Since $(z_\ep,\xi_\eps)_{\eps\in (0,1]}$ is 
bounded, there exists a sequence $(\eps_n)_n$ such that
$$
 \eps_n\downarrow 0
 \quad \mbox{and}\quad 
 (z_n,\xi_n):=\big(z_{\eps_n},\xi_{\eps_n}\big)
 \longrightarrow (\hat z,\hat\xi)=(z_0,\hat\xi) \in(0,\infty)\times\R,
$$
where the fact that $\hat z=z_0$ follows from the strict maximum 
property in \reff{strictmax} and  classical arguments from the 
theory of viscosity solutions.
We finally conclude from \reff{subsol-almostfinal} and 
\reff{remainderestimate-subsol} that
 \b*
 0 
 &\ge&
 - \frac12(\sigma^2 v_{zz})(s_0,z_0)\hat\xi^2
 -\Ac\varphi(s_0,z_0)
 -\Ac\phi(0)
 +\frac12(1-\delta)\alpha^2(s_0,z_0)w_{\xi\xi}(s_0,z_0,\hat\xi)
 \\
 &=&
 -\Ac\varphi(s_0,z_0)
 - \frac12(\sigma^2 v_{zz})(s_0,z_0)\hat\xi^2
 +\frac12(1-\delta)\alpha^2(s_0,z_0)w_{\xi\xi}(s_0,z_0,\hat\xi),
 \e*
since $\Ac\phi(0)=0$. Now, in view of the first corrector 
equation \reff{e.cw}, 
 \begin{eqnarray*}
 0 
 &\ge&
 -\Ac\varphi(s_0,z_0)
 +a(s_0,z_0)
 +\frac12\delta\alpha^2(s_0,z_0)w_{\xi\xi}(s_0,z_0,\hat\xi).
 \e*
Finally, we conclude that $\Ac\varphi(s_0,z_0) - a(s_0,z_0)\ge 0$, 
by sending $\delta$ to zero.
\qed

\section{Verifying Assumption \ref{a.bound}}
\label{s.verify}

In this section, we verify Assumption 
\reff{a.bound}.
This is done
by constructing
an appropriate sub-solution
of the dynamic programming equation \reff{e.dpp}.
Clearly, this construction requires
assumptions and here we present only one
possible set of assumptions.
To simplify the presentation,
we suppose that the coefficients are independent
of the $s$-variable.  
Next, we assume that
there exist constants $0< k_*\le k^*$ so that 
the limit Merton value function
satisfies
\begin{eqnarray}
\label{a.almost}
0< k_*z \le \eta(z) \le k^* z.
\end{eqnarray}
Let $ \cbf$ be the optimal Merton consumption policy
given as in \reff{hatc}.  We assume that 
\be\label{a.cpower}
U(\cbf(z)) \ge k_* z\vp(z),
\ee
for some constant $k_* >0$.  Notice that
all the above assumptions hold in the power utility case.
First, using \reff{e.alphabound} and
the explicit representation of $a$, 
one may directly verify that
there is a constant $a^*>0$
so that
$$
a(z) \le a^* z \vp(z).
$$
Then, the definition of $\Ac$ and the above assumptions 
imply that
\begin{equation}
\label{e.abound}
\Ac v(z) = U(\cbf(z)) \ge k_* z\vp(z) \ge 
\frac{k_*}{a^*} a(z) =   \frac{k_*}{a^*} \Ac u(z).
\end{equation}
Let $u$ be the function defined in \reff{e.defu}. Since $v$ is assumed to be smooth, we may apply It\^o's formula in a standard way to conclude from the last inequality that
\begin{equation}
\label{e.vu}
0 \le u(z) \le \frac{a^*}{k_*}\ v(z) .
\end{equation}
Moreover, since we assume that coefficients
are independent of the $s$ variable, 
\reff{hattheta} is equivalent to $\ybf(z) =\eta(z) (\mu-r)/\sigma^2$.
Hence, \reff{e.alphabound}  implies
that
\begin{equation}
\label{e.third}
- \vpp(z) \le \eta(z) \ v^{\prime \prime \prime} \le - 2\vpp(z).
\end{equation}

We now use these observations 
to construct a sub-solution of
the dynamic programming equation of the 
form
\begin{equation}
\label{e.ve}
V^\ep(x,y):= v(z) -K \ep^2 v(z) + \ep^4 
\tilde{W}(z,\xi),
\end{equation}
with a sufficiently large constant $ K \ge a^*/k_*$
and a slightly modified corrector,
$$
\tilde{W}(z,\xi):= z \vp(z) \tilde{w}(\xi/z),
$$ 
where the function $\tilde{w}(z)$
and the constant $\tilde{a}>0$ are the unique
solution of $\tilde{w}(0)=0$ and
\begin{equation}
\label{e.tildew}
\max \left\{
-\frac{k_* \sigma^2}{2} \rho^2 - \frac{(\alpha^* k^* )^2}{2}
\tilde{w}_{\rho \rho} +\tilde{a}\ ;\
-2 \lambda^{1,0}+ \tilde{w}_\rho\ ;\ 
-2 \lambda^{0,1}- \tilde{w}_\rho \right\} .
\end{equation}
The solution of the above equation is explicitly available
through the general solution obtained earlier in Section \ref{s.w1d}.

The fact that $V^\ep$ is a sub-solution of \reff{e.dpp}
follows from  tedious but otherwise direct calculations.
To streamline these calculations, we first state an 
estimate that follows from the explicit form of $\tilde{W}$.

\begin{Lemma}
\label{l.estimates}
There is a constant $k^*>0$ so that
\begin{eqnarray*}
&z \left|\tilde{W}_{\xi \xi}(z,\xi)\right| \le
k^*  \vp(z), \\
& \left|\tilde{W}_z(z,\xi) \right| \le 
k^*  \vp(z) \left(1+ \frac{|\xi|}{z}\right), \\
& z\left|\partial_x \tilde{W}(z,\xi) \right|
+z\left| \partial_y  \tilde{W}(z,\xi) \right|  \le 
k^* z \vp(z) \left(\frac{1}{\ep}+ \frac{|\xi|}{z}\right), \\
& z^2 \left| \partial_{y y}\tilde{W}(z,\xi) 
-\frac{(1-\ybf^\prime(z))^2}{\ep^2}  
\tilde{W}_{\xi \xi }(z,\xi) \right|  \le 
k^* z \vp(z) \left(\frac{1}{\ep}+ \frac{|\xi|}{z}\right).
\end{eqnarray*}
\end{Lemma}

\no
{\em{Proof.}} These estimates follow directly from straightforward
differentiation and the estimates \reff{a.almost},
\reff{e.third}.
\qed

\begin{Lemma}[Lower Bound]
\label{l.upper}
Assume \reff{a.almost}, \reff{a.cpower} and \reff{a.theta}.  
Then, for sufficiently large $K>0$,
$V^\ep$ defined in \reff{e.ve} is a 
sub-solution of \reff{e.dpp}
in $\R^2_+$.  Moreover,
$$
\bar{u}^\ep(x,y) \le K v(z) + \ep^2 \tilde{W}(z,\xi)
$$
on $\R^2_+$
and Assumption
\ref{a.bound} holds.
\end{Lemma}

\no
{\em{Proof.}}
We need to show that 
at any point $(x,y)\in \R^2_+$
one of the three terms in \reff{e.dpp} is non-positive.
Since $(x,y)\in \R^2_+$, 
by assumption \reff{a.theta}, we have
$$
|\xi|= \frac{|y-\ybf(z)|}{\ep} \le \frac{z}{\ep},
\quad
\Rightarrow
\quad
\Xi := \frac{\xi}{z} \in \frac{1}{\ep}\ [-1,1].
$$
Let $\rho_0>0$ be the threshold in the 
equation \reff{e.tildew}.
We analyze several cases separately. 

\no
{\em{Case 1.} } $\rho_0 \le \Xi \le 1/\ep$.

In this case,
$\tilde{W}_\xi(z,\xi)=2\lambda^{1,0}\vp(z)$.
We 
use the previous Lemma and \reff{a.theta},
to arrive at,
\begin{eqnarray*}
\Lambda^\ep_{1,0}\cdot (V^\ep_x,V^\ep_y)
&=&  \frac{1}{\ep}\hat{V}^\ep_\xi 
+ \ep^2 \lambda^{1,0}  (1- \ybf^\prime)\hat{V}^\ep_\xi
+ \ep^3\lambda^{1,0} \hat{V}^\ep_z\\
&=& \ep^3 \left[(1-\ep^3 \lambda^{1,0} (1-\ybf^\prime)) 
\tilde{W}_\xi 
+ (1-C\ep^2 ) \vp  -\lambda^{1,0}\ep^4 
\tilde{W}_z\right]\\
&\le & \ep^3 \lambda^{1,0}\vp \left(-1 +k^*\ep^3
\right)  \le 0,
\end{eqnarray*}
provided that
$\ep $ is sufficiently small.
\vspace{10pt}

\no
{\em{Case 2.} }  $- 1/\ep\le \Xi \le -\rho_0 $.

A similar calculation,
shows that $\Lambda^\ep_{0,1}\cdot (V^\ep_x,V^\ep_y) \le 0$, for all
sufficiently small $\ep$.
\vspace{10pt}

\no
{\em{Case 3.} }  $|\Xi| \le \rho_0$.
We now use Remark \ref{r.formal} to conclude that
$$
\Jc(V^\ep) = \ep^2 \left[   -\frac{\sigma^2\vpp(z)}{2} \xi^2 
 +\frac{\alpha^2(z)}{2}\tilde{W}_{\xi \xi}(z,\xi) -
  K \Ac v(z)+\Rc^\ep(z,\xi)\right].
$$
We first use  \reff{a.almost}, \reff{a.theta},
\reff{e.tildew}, \reff{e.abound}
and set $\rho:= \xi/z$.
The result is
\begin{eqnarray*}
\Ic &:= &\frac{\Jc(V^\ep)}{\ep^2}\\
 &\le & \ep^2 \vp(z)\eta(z)
\left[  \frac{k_*\sigma^2}{2} \rho^2 
 +\frac{(\alpha^* k^*)^2}{2}\tilde{w}_{\rho \rho}(\rho) 
 -K(k_*)^2\right] +\ep^2 \Rc^\ep(z,\xi)\\
&= & \ep^2 \vp(z)\eta(z)
\left[  \tilde{a}
 -K(k_*)^2\right] +\ep^2 \Rc^\ep(z,\xi).
\end{eqnarray*}
If $K$ is sufficiently large then $K(k_*)^2$ is larger than $\tilde{a}$
and by \reff{a.almost}, the above estimate implies that
$$
\Ic  \le - z \vp(z) + \Rc^\ep(z,\xi).
$$ 
We now estimate $\Rc^\ep$ by recalling the results of
subsection  \ref{s.remainder}.  We split this in three terms
coming from the value function $v$, the corrector $\tilde{W}$
and from the utility function,
$$
|\Rc^\ep|:= \Rc^\ep_v + \Rc^\ep_w + \Rc^\ep_U.
$$
We estimate each one using Lemma \ref{l.estimates}. Then,
\begin{eqnarray*}
\Rc^\ep_v& \le &K \left[ \ep\Xi (\mu -r) z \vp(z)+ \frac{\sigma^2}{2} 
\left( \ep^2\Xi^2 +  2\ep \Xi  (\ybf/z)\right) z^2 \vpp(z)\right]\\
& \le &\ep K k^* z \vp(z).
\end{eqnarray*}
Also
\begin{eqnarray*}
\Rc^\ep_w& \le & \ep^2\Big[\beta \tilde{W} 
- rz \left((1-(\ybf/z)) + \ep \Xi \right)\tilde{W}_x
-\mu z \left(\ep\Xi+ (\ybf/z)\right) \tilde{W}_y \\
&& - \frac{\sigma^2}{2} z^2 \left(\ep\Xi+ (\ybf/z)\right) ^2
\left(\tilde{W}_{yy} - \tilde{W}_{\xi \xi} (1-\ybf_z)^2/\ep^2\right) \\
&&  \frac{\sigma^2}{2} z^2
\tilde{W}_{\xi \xi}\frac{ (1-\ybf_z)^2}{\ep^2}
 \left(\ep^2\Xi^2+ 2\ep \Xi (\ybf/z)\right) 
\\
& \le &k^* z \vp(z).
\end{eqnarray*} 
Finally
\begin{eqnarray*}
\Rc^\ep_U &=& \tilde{U}(\vp) - 
\tilde{U}(V^\ep_x)\\
& \le & 
\tilde{U}(\vp) - 
\tilde{U}(\vp[1-\ep^2K +k^* \ep^4]) \le 0.
\end{eqnarray*} 
Hence, there is $k^*$ so that.
 $$
| \Rc^\ep| \le \ep k^*z \vp(z) .
$$ 
Hence if $K$ is sufficiently large, $V^\ep$ is a sub-solution 
of \reff{e.dpp} for all small $\ep$.
\vspace{10pt}

\no
{\em {Boundary}} $y=0$.

Then, again by 
\reff{a.theta}, for all sufficiently small $\ep>0$,
$$
\Xi = \frac{y -\ybf(z)}{\ep} =\frac{-\ybf(z)}{\ep}
 <-\rho_0.
$$
Hence, by the second case, and  Lemma \ref{l.notrade}
$$
\Lambda^\ep_{1,0}\cdot (V^\ep_x,V^\ep_y)(x,0) \le 0 = 
\Lambda^\ep_{1,0}\cdot (v^\ep_x,v^\ep_y)(x,0),
\quad
\forall\ x >0.
$$  

\no
{\em {Boundary}} $x=0$.

By a similar analysis, we can show that
$$
\Lambda^\ep_{0,1}\cdot (V^\ep_x,V^\ep_y) (0,y) \le 0 = 
\Lambda^\ep_{0,1}\cdot (v^\ep_x,v^\ep_y)(0,y),
\quad
\forall\ y >0.
$$  
Then, on $\R^2_+$, $V^\ep$ is a sub-solution of \reff{e.dpp}
while $\ve$ is a solution.  Also on the boundary of $\R^2_+$
again $V^\ep$ is a sub-solution of an oblique Neumann condition
and $\ve$ is a super-solution.  Then, by comparison (or by a 
verification argument), we conclude that
$\ve \ge \phi$ on $\R^2_+$.  This proves the lower bound on $\ue$
on the positive orthant. 
\qed

\begin{Remark}\label{rem:behaviorzero}{\rm
In view of Lemma \ref{l.upper}, it follows that the local upper bounding function $B$, defined in \reff{e.B}, is bounded by the function $K v(z)$. In particular, this implies that the growth of $u_*$ and $u^*$, both at infinity and at the origin, is the same as that of the zero-transaction cost Merton value function $v$. By introducing the logarithmic variable, we observe that the behavior near the origin transforms into a growth condition at minus infinity.
}
\end{Remark}

\section{Homothetic case}
\label{s.homothetic}
In this short section, we consider the classical CRRA utility function
 \be\label{powerutility}
U(c) := \frac{c^{1-\gamma}}{1-\gamma},
\qquad c >0,
 \ee
for some $\gamma >0$ with $\gamma =1$ corresponding 
to the logarithmic utility. Our objective is to reproduce the results of 
Janecek and Shreve \cite{js2004} by directly applying our explicit 
expansion result of Theorem \ref{t.main}. Also these
calculations show how one may use our
results to obtain the asymptotic formulae
for problems with power utility that
have explicitly known Merton value functions,
such as factor models.

In the context of the power utility \reff{powerutility}, 
the Merton value function is explicitly given by,
$$
v(z)= \frac{1}{(1-\gamma)} \  \frac{z^{1-\gamma}}{v_M^{\gamma}}, 
$$
with the Merton constant
$$
v_M= \frac{\beta-r(1-\gamma)}{\gamma}
- \frac12 \frac{(\mu-r)^2}{\gamma^2 \sigma^2}(1-\gamma).
$$
Hence the risk tolerance function and 
the optimal strategies are given by,
$$
\eta(z)= \frac{z}{\gamma},
\qquad
\ybf(z)= \frac{\mu-r}{\gamma \sigma^2}\ z:= \pi_M z,
\qquad
\cbf(z) = v_M z.
$$
In particular, since $\ybf$ and $\cbf$ are linear in $z$, the comparison 
Assumption \ref{a.compare} is immediately checked to hold true. Indeed, by introducing the logarithmic variable $z'=\ln{z}$, the second corrector equation \reff{e.au} becomes linear with constant coefficients on $(-\infty,\infty)$. The growth condition as discussed in Remark \ref{rem:behaviorzero} transforms into an exponential sublinear growth. It is well-known that this condition is sufficient to prove comparison. The corresponding probabilistic argument refers to the integrability of exponential sublinear growth with respect to the Gaussian density.

Moreover, since the conditions of Section \ref{s.verify} are satisfied in 
the present context, it follows that Assumptions \ref{a.bound} holds 
true in our power utility case, provided that $\pi_M\in(0,1)$. Finally, by 
Remark 11.3 in Shreve and Soner \cite{ss}, the last condition also 
implies the validity of Assumption \ref{a.merton}. 
We have then verified the following.

\begin{Lemma}
Assume $\pi_M\in(0,1)$. Then, Assumptions \ref{a.bound}, \ref{a.smooth},
 \ref{a.compare} and \ref{a.merton} hold true in the context of the power 
 utility function \reff{powerutility}. 
\end{Lemma}
 
Since the diffusion coefficient $\alpha(z)=\sigma \ybf(z)[1-\ybf_z(z)]$, it follows that
$$
\bar{\alpha}= \frac{\alpha(z)}{\eta(z)}=\gamma \sigma \pi_M (1-\pi_M).
$$
The constants in the solution of the corrector equation are given by,
$$
\rho_0= \left(\frac{3\bar{\alpha}^2}{4 \sigma^2}\left(\lambda^{1,0}+
\lambda^{0,1}\right)\right)^{1/3},
$$
$$
a(z) = \eta(z)\vp(z) \bar{a}= \frac{\sigma^2(1-\gamma)}{2 \gamma} 
 \rho_0^2 \ v(z). 
$$
Since
$$
\Ac v(z) = U(\cbf(z)) = \frac{1}{1-\gamma} \left( v_M z\right)^{1-\gamma}
= v_M v(z),
$$
the unique solution $u(z)$ of the second corrector equation 
$$
\Ac u(z) =a(z) =  \frac{\sigma^2(1-\gamma)}{2 \gamma} 
 \rho_0^2 \ v(z)
$$
is given by
$$
u(z) = \frac{\sigma^2(1-\gamma)}{2 \gamma} 
 \rho_0^2 v_M^{-1}\ v(z)
= u_0 z^{1-\gamma},
$$
where
$$
u_0:= 
\left( \pi_M(1-\pi_M)\right)^{4/3} v_M^{-(1+\gamma)} .
$$
Finally, we summarize the  expansion result in the following.

\begin{Lemma}
\label{l.homothetic}
For the power utility function $U$ in \reff{powerutility}, 
$$
v^\ep(x,y) = v(z) - \ep^2 u_0 z^{1-\gamma}+O(\ep^3).
$$
The width of the transaction region
for the first correction equation $2\xi_0=2\eta(z)\rho_0$
is given by
$$
2\xi_0= \left( \frac{6}{\gamma}(\lambda^{0,1}+\lambda^{1,0})\right)^{1/3}   
\left( \pi_M (1-\pi_M)\right)^{2/3}.
$$
\end{Lemma}
The above formulae with $\lambda^{i,j}=1$
are exactly the same 
as equation (3.13) in Janecek and Shreve \cite{js2004} .

\begin{paragraph}{Acknowledgements}

\end{paragraph}

\end{document}